\documentclass[final,leqno,letterpaper]{etna}
\usepackage{siunitx}  
\usepackage{upgreek}
\usepackage[intlimits]{amsmath}
\usepackage{amsfonts}
\usepackage{graphicx}
\usepackage[linesnumbered,ruled,vlined]{algorithm2e}
\usepackage{flafter}
\usepackage{booktabs}
\usepackage{threeparttable}
\usepackage{multirow}
\usepackage{multicol}
\usepackage{tikz} 
\usepackage{float}
\usepackage[space,noadjust,nocompress]{cite}

\setbibdata{1}{xx}{xx}{2022} 

\hypersetup{%
    pdftitle={Original Manuscript},
    pdfauthor={Jinqiang Chen, Vandana Dwarka, Cornelis Vuik},
    pdfkeywords={Helmholtz equation, parallel computation, and matrix-free}
    }

\title{A matrix-free parallel solution method for the three-dimensional heterogeneous Helmholtz equation\thanks{%
Received... Accepted... Published online on... 
}}

\author{Jinqiang Chen\footnotemark[2]
        \and Vandana Dwarka\footnotemark[2]
        \and Cornelis Vuik\footnotemark[2]}

\shorttitle{MARIX-FREE PARALLEL METHOD FOR 3D HELMHOLTZ} 
\shortauthor{J. CHEN, V. DWARKA AND C.VUIK}

\begin{document}

\maketitle

\renewcommand{\thefootnote}{\fnsymbol{footnote}}

\footnotetext[2]{Department of Applied Mathematics, Delft University of Technology, Delft, the Netherlands
(\tt{\{j.chen-11, v.n.s.r.dwarka, c.vuik\}@tudelft.nl})}

\begin{abstract}
The Helmholtz equation is related to seismic exploration, sonar, antennas, and medical imaging applications. It is one of the most challenging problems to solve in terms of accuracy and convergence due to the scalability issues of the numerical solvers. For 3D large-scale applications, high-performance parallel solvers are also needed. In this paper, a matrix-free parallel iterative solver is presented for the three-dimensional (3D) heterogeneous Helmholtz equation. We consider the preconditioned Krylov subspace methods for solving the linear system obtained from finite-difference discretization. The Complex Shifted Laplace Preconditioner (CSLP) is employed since it results in a linear increase in the number of iterations as a function of the wavenumber. The preconditioner is approximately inverted using one parallel 3D multigrid cycle. For parallel computing, the global domain is partitioned blockwise. The matrix-vector multiplication and preconditioning operator are implemented in a matrix-free way instead of constructing large, memory-consuming coefficient matrices. Numerical experiments of 3D model problems demonstrate the robustness and outstanding strong scaling of our matrix-free parallel solution method. Moreover, the weak parallel scalability indicates our approach is suitable for realistic 3D heterogeneous Helmholtz problems with minimized pollution error. 
\end{abstract}



\section{Introduction} 

    The Helmholtz equation, describing the phenomena of time-harmonic wave scattering in the frequency domain, finds applications in many scientific fields, such as seismic problems, advanced sonar devices, medical imaging, and many more. To solve the Helmholtz equation numerically, we discretize it and obtain a linear system $Ax=b$. The linear system matrix is sparse, symmetric, complex-valued, non-Hermitian, and indefinite. Instead of a direct solver, iterative methods and parallel computing are commonly considered for a large-scale linear system resulting from a 3D problem. However, the indefiniteness of the linear system brings a great challenge to the numerical solution method, especially for large wavenumbers. The convergence rate of many iterative solvers is affected significantly by increasing wavenumber. An increase in the wavenumber leads to a dramatic increase in iterations. Moreover, the general remedy for minimizing the so-called pollution error, driven by numerical dispersion errors due to discrepancies between the exact and numerical wavenumber, is to refine the grid such that the condition $k^3h^2 < 1$ is satisfied \cite{babuska1997pollution}. Therefore, the research problem of how to solve the systems efficiently and economically while at the same time maintaining a high accuracy by minimizing pollution error arises in this field. A wavenumber-independent-convergent and parallel scalable iterative method could significantly enhance the corresponding research in electromagnetics, seismology, and acoustics.
    
    Many efforts have been made to solve the problem in terms of accuracy and scalable convergence behavior. One of the main concerns is the spectrum of the system matrix, which is closely related to the convergence of Krylov subspace methods. The preferable idea is to preprocess the system with a preconditioner. By applying a preconditioner to the linear system, the solution remains the same, but the coefficient matrix has a more favorable distribution of eigenvalues. 
	Many preconditioners are proposed for the Helmholtz problem so far, such as incomplete Cholesky (IC)/incomplete LU (ILU) factorization, shifted Laplacian preconditioners \cite{bayliss1983iterative,erlangga2004class,erlangga2006novel} and so on. The industry standard, also known as the Complex Shifted Laplace Preconditioner (CSLP) \cite{erlangga2004class, erlangga2006novel} does show good properties for medium wavenumbers. Nevertheless, the eigenvalues shift to the origin as the wavenumber increases. These near-zero eigenvalues have an unfavorable effect on the convergence speed of the Krylov-based iterative solvers. Therefore, Sheihk \cite{sheikh2016accelerating} further included deflation techniques to accelerate the convergence. Recently, a higher-order approximation scheme to construct the deflation vectors was proposed to reach close to wavenumber-independent convergence \cite{dwarka2020scalable}. 

	The development of scalable parallel Helmholtz solvers is also ongoing. One approach is to parallelize existing advanced algorithms. Kononov and Riyanti \cite{Kononov2007Numerical,riyanti2007parallel} first developed a parallel version of Bi-CGSTAB preconditioned by multigrid-based CSLP.  Knibbe et al. \cite{knibbe2011gpu} further introduced parallel versions of CSLP-preconditioned Bi-CGSTAB and IDR(s) which run on GPU accelerators. Dan and Rachel \cite{Gordon2013Robust} parallelized their so-called CARP-CG algorithm (Conjugate Gradient acceleration of CARP) blockwise. The block-parallel CARP-CG algorithm shows improved scalability as the wavenumber increases. Calandra et al. \cite{calandra2013improved,calandra2017geometric} proposed a geometric two-grid preconditioner for 3D Helmholtz problems, which shows the strong scaling property in a massively parallel setup.

	Another approach is the Domain Decomposition Method (DDM), which originates from the early Schwarz Methods. DDM has been widely used to develop parallel solution methods for Helmholtz problems. For comprehensive surveys, we refer the reader to \cite{douglas1998second,mcinnes1998additive,toselli1999overlapping,collino2000domain,gander2002optimized,schadle2007additive,engquist2011sweeping,boubendir2012quasi,chen2013source,stolk2013rapidly,gander2019class,taus2020sweeps} and references therein.
	
	This work is interested in parallelizing Krylov subspace methods, such as the Generalized Minimal RESidual method (GMRES), Bi-CGSTAB, and IDR(s), preconditioned by the multigrid-based CSLP for the Helmholtz equation. We consider the CSLP preconditioner because it is the first and most popular method where the number of iterations scales linearly within medium wavenumbers. However, CSLP is not wavenumber independent. The standard configuration of CSLP+GMRES works less efficiently if we want to solve highly heterogeneous problems with minimized pollution error. Finer grids and more iterations, hence more memory, will be needed. Our contribution is the development and validation of a matrix-free parallel framework of CSLP-preconditioned Krylov subspace methods in the context of solving large-scale 3D Helmholtz problems with minimized pollution error. To the best of our knowledge, this has not been previously reported in the literature. The earlier variants proposed by Kononov and Riyanti et al. \cite{Kononov2007Numerical,riyanti2007parallel} mainly parallelized the sequential program based on the data-parallel concept. It results in a row-wise domain partition and a 3D multigrid method with 2D semi-coarsening. In contrast, this work starts with a block-wise domain partition and implements a standard 3D multigrid method in a matrix-free way. Our method contributes to a robust and scalable parallel CSLP-preconditioned solver for realistic 3D applications. Numerical experiments on typical 3D model problems show that the matrix-free parallel solution method can effectively save memory for storing global sparse matrices and show good parallel performance. 
	
	The rest of this paper is organized as follows. Section \ref{sec:mathmodels} describes the mathematical model and discretization technique that we will discuss. All numerical methods we use are described in Section \ref{sec:methods}. Section \ref{sec:implem} describes the parallel implementation. The numerical performance is explored in Section \ref{sec:exper}. Finally, Section \ref{sec:concls} contains our conclusions.

\section{Mathematical models} \label{sec:mathmodels}
    The Helmholtz equation can model the wavefield in heterogeneous media in the frequency domain. Suppose a parallelepipedal domain $\Omega \subset \mathbb{R}^3 $ with boundary $\Gamma = \partial\Omega$. The Helmholtz equation reads
	
    \begin{equation}\label{eq:HelmholtzEq}
    -\Delta u(\mathbf{x}) - k(\mathbf{x})^2 u(\mathbf{x}) = b(\mathbf{x}),\  \mathbf{x}=(x_1,x_2,x_3) \in \Omega, \nonumber
    \end{equation}
    supplied with one of the following boundary conditions:
    \begin{equation}\label{eq:DirichletBC}
    \text{Dirichlet: } u(\mathbf{x})=g(\mathbf{x}),\  \text{on} \ \partial\Omega,  \nonumber
    \end{equation}
	
    \begin{equation}\label{eq:SommerfeldBC}
    \text{First-order Sommerfeld: } \frac{\partial u(\mathbf{x})}{\partial \vec{n}}-\text{i} k(\mathbf{x}) u(\mathbf{x}) = \mathbf{0},\  \text{on} \ \partial\Omega, 
    \end{equation}
    where $\text{i}$ is the imaginary unit. $u$ represents the pressure wavefield. $b$ is the source function. $\vec{n}$ and $g$ represent the outward unit normal vector and the given data on the boundary, respectively. $k(\mathbf{x})$ is the wavenumber, suppose the frequency is $f$, they are related by
    \begin{equation}\label{eq:wavenum}
    k(\mathbf{x}) =\frac{2\pi f}{c(\mathbf{x})}. \nonumber
    \end{equation}
    where $c(\mathbf{x})$ is the space-variant acoustic-wave velocity due to changes in the material of the domain.
    
    \subsection{3D closed-off problem}
    This problem is a constant wavenumber model with a given solution to validate the numerical methods. A parallelepipedal homogeneous domain $\Omega=\left[0,1 \right]^3$ is considered. The source function is specified by
    \begin{equation}\label{eq:3Dcloseoff}
    b\left( {x_1,x_2,x_3} \right) = \left( {21{\pi ^2} - {k^2}} \right)\sin \left( {\pi x_1} \right)\sin \left( {2\pi x_2} \right)\sin \left( {4\pi x_3} \right) - {k^2},\ \Omega=\left[0,1 \right]^3.  \nonumber
    \end{equation}
    It is supplied with the following Dirichlet conditions
    \begin{equation}
    u(x_1,x_2,x_3)=1,\  \text{on} \ \partial\Omega. \nonumber
    \end{equation}
    The analytical solution is given by
    \begin{equation}
    \label{eq: MP1 analytical solution}
    u(x_1,x_2,x_3) = \sin \left( {\pi x_1} \right)\sin \left( {2\pi x_2} \right) \sin(4 \pi x_3) + 1. \nonumber
    \end{equation}
    
    \subsection{3D Wedge problem} \label{sec:wedge}
\begin{figure}[htbp]
        \centering
        \includegraphics[width=0.6\textwidth,trim=2 2 2 2,clip]{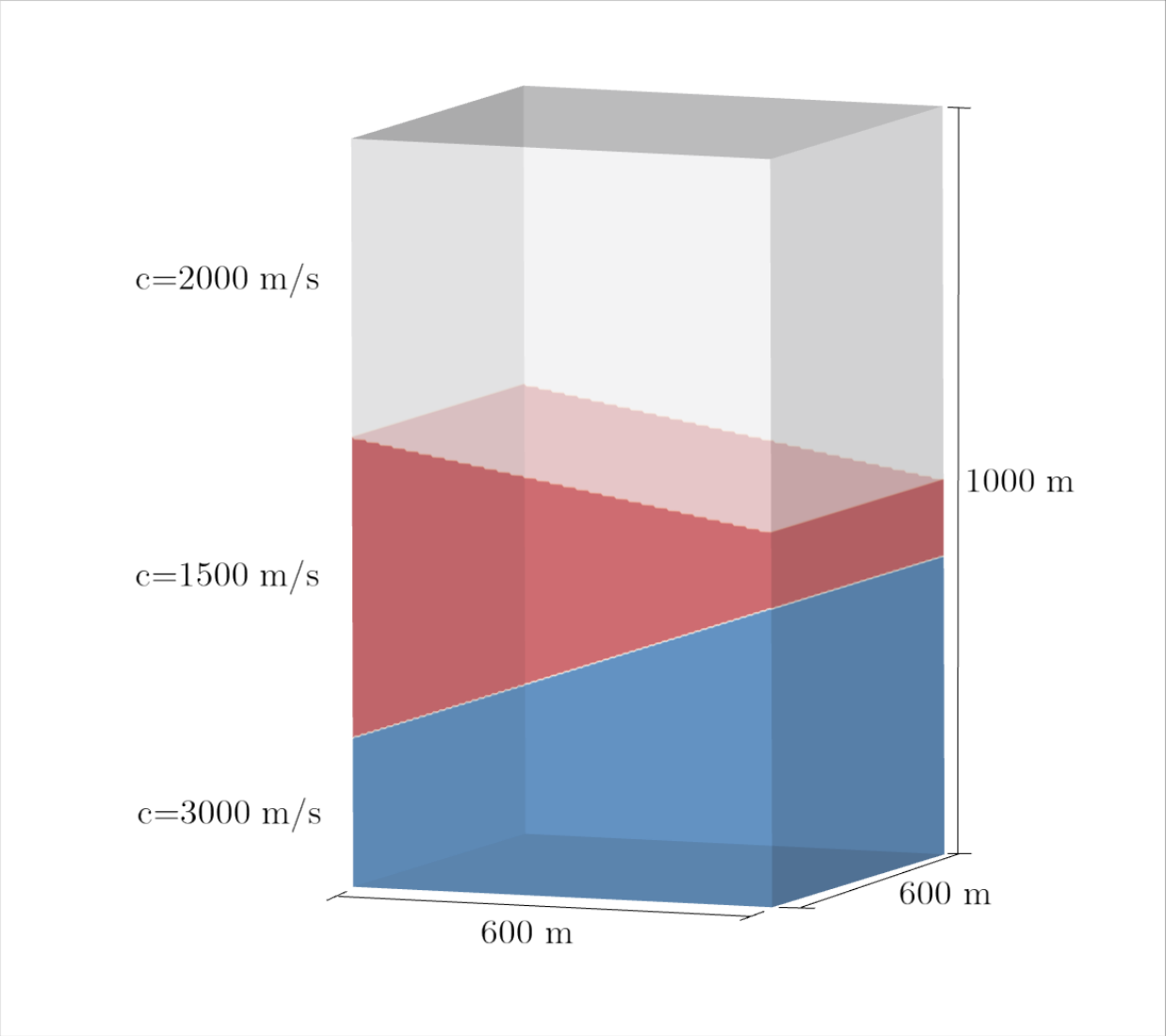}
        \caption{The velocity distribution of the 3D wedge problem.}
        \label{fig:3D wedge domain}
    \end{figure}

    Most physical problems of geophysical seismic imaging describe a heterogeneous medium. The so-called Wedge problem \cite{plessix2003separation} is a typical problem with a simple heterogeneity. It mimics three layers with different acoustic-wave velocities hence different wavenumbers. As shown in Figure \ref{fig:3D wedge domain}, the parallelepipedal domain $\Omega=\left[0,600 \right] \times \left[0,600 \right] \times \left[-1000,0 \right]$  is split into three layers. Suppose the acoustic-wave velocity $c$ is constant within each layer but different from each other. A point source is located at $\left(x_1, x_2, x_3 \right) = \left( 300, 300, 0\right)$. 
    
    The problem is given by
    \begin{equation}\label{eq: 3D wedge problem}
        \left\{ {\begin{array}{*{20}{c}}
            {-\Delta u(x_1,x_2,x_3)-k(x_1,x_2,x_3)^2u(x_1,x_2,x_3)=b(x_1,x_2,x_3), \quad\text{on } \Omega}\\
            {b(x_1,x_2,x_3)=\delta(x_1-300,x_2-300,x_3) \quad x_1,x_2,x_3\in \Omega}
            \end{array}} \right. , \nonumber
    \end{equation}
    where $k(x_1,x_2,x_3)=\frac{2\pi f}{c(x_1,x_2,x_3)}$. $f$ is the frequency. The wave velocity $c(x_1,x_2,x_3)$ is shown in Figure \ref{fig:3D wedge domain}. $\delta \left( {x_1,x_2,x_3} \right)$ is a Dirac delta function. The first-order Sommerfeld boundary conditions are imposed on all boundaries. 
    
    \subsection{3D SEG/EAGE salt model}
\begin{figure}[htbp]
        \centering
        \includegraphics[width=0.85\textwidth,trim=2 2 2 2,clip]{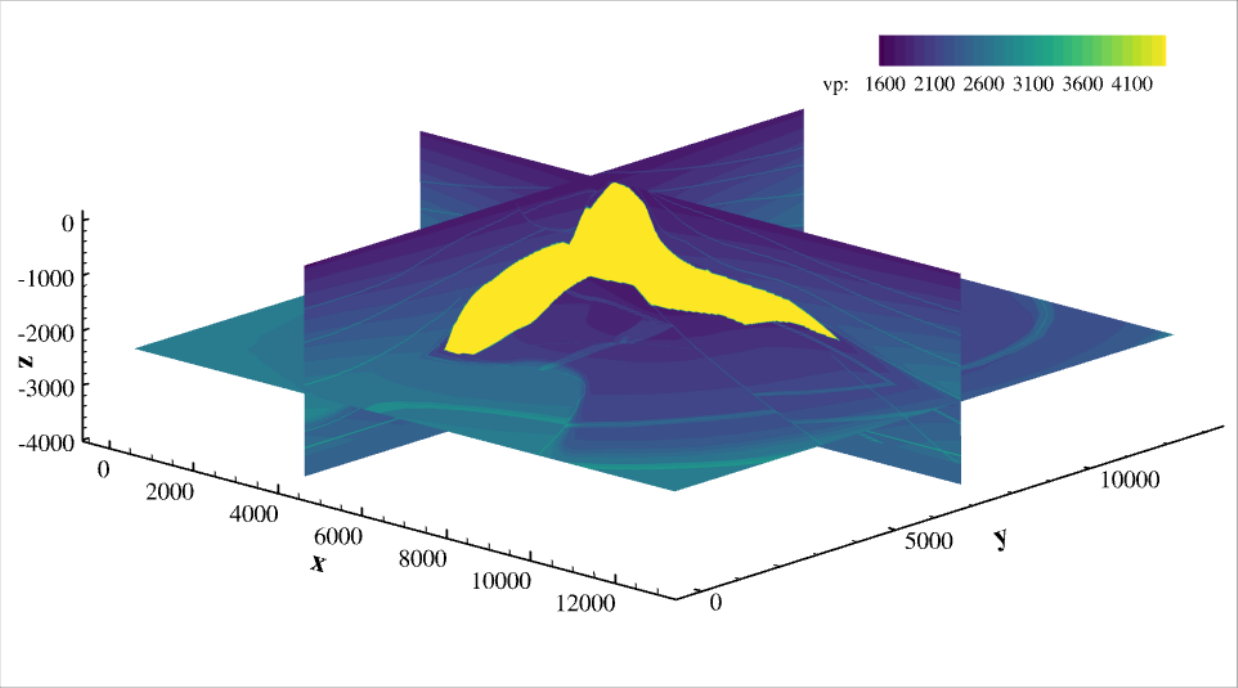}
        \caption{The velocity distribution of the 3D SEG/EAGE Salt Model. The velocity varies from $\SI{1500}{\metre\per\second}$ to $\SI{4482}{\metre\per\second}$.}
        \label{fig:3D Salt Model}
    \end{figure}

    The 3D SEG/EAGE salt model \cite{aminzadeh1997seg} is a velocity field model containing salt domes, which mimics the typical Gulf Coast salt structure. As shown in Figure \ref{fig:3D Salt Model}, it is defined in a parallelepipedal physical domain of size $ \SI{13520}{\metre} \times \SI{13520}{\metre} \times \SI{4200}{\metre}$. The acoustic-wave velocity varies from $\SI{1500}{\metre\per\second}$ to $\SI{4482}{\metre\per\second}$. The model is considered challenging due to the inclusion of complex geometries (salt domes) and a realistic large-size computational domain. We consider a computational domain $\Omega=(0,12800)\times(0,12800)\times(0,3840)$ with grid size $641\times641\times193$, allowing multi-level geometric coarsening. A point source is located at $\left(x_1, x_2, x_3 \right) = \left( 3200, 3200, 0\right)$. 
    
    The problem is given by
    \begin{equation}\label{eq: 3D Salt problem}
        \left\{ {\begin{array}{*{20}{c}}
            {-\Delta u(x_1,x_2,x_3)-k(x_1,x_2,x_3)^2u(x_1,x_2,x_3)=b(x_1,x_2,x_3), \quad\text{on } \Omega}\\
            {b(x_1,x_2,x_3)=\delta(x_1-3200,x_2-3200,x_3) \quad x_1,x_2,x_3\in \Omega}
            \end{array}} \right. , \nonumber
    \end{equation}
    where $k(x_1,x_2,x_3)=\frac{2\pi f}{c(x_1,x_2,x_3)}$. $f$ is the frequency. The wave velocity $c(x_1,x_2,x_3)$ is shown in Figure \ref{fig:3D Salt Model}. The first-order Sommerfeld boundary conditions are imposed on all boundaries. 
    
    \subsection{Finite-difference discretization}
    Structural vertex-centered grids are used to discretize the computational domains. Suppose the grid widths in $x_1$, $x_2$, and $x_3$ directions are all equal to $h$.  A second-order finite difference scheme for a 3D Laplacian operator has the following stencil
    \begin{equation}\label{eq:stencilLaplace}
    \left[ {{-\Delta_h}} \right] = \frac{1}{{{h^2}}}
    \left[ 
        \left[\begin{array}{*{20}{c}}
            0&0 &0\\
            0&-1&0\\
            0&0 &0
            \end{array}\right]_{x_3-h}
        \left[\begin{array}{*{20}{c}}
        0&{ - 1}&0\\
        { - 1}&{6}&{ - 1}\\
        0&{ - 1}&0
        \end{array}\right]_{x_3}
        \left[\begin{array}{*{20}{c}}
            0&0 &0\\
            0&-1&0\\
            0&0 &0
        \end{array}\right]_{x_3+h} 
    \right].
    \end{equation}
    The discrete Helmholtz operator $A_h$ can be obtained by adding the diagonal matrix $-k^2 I_h$  to the Laplacian operator $-\Delta_h$, i.e.,
    \begin{equation}\label{eq:operator}
    A_h =-\Delta_h-k^2 I_h. \nonumber
    \end{equation}
    Therefore, the stencil of the discrete Helmholtz operator is
    \begin{equation}\label{eq:Helmstencil}
    \left[ {{A_h}} \right] = \frac{1}{{{h^2}}}\left[
    \left[ {\begin{array}{*{20}{c}}
        0&0 &0\\
        0&-1&0\\
        0&0 &0
    \end{array}} \right]_{x_3-h}	
    \left[ {\begin{array}{*{20}{c}}
        0&{ - 1}&0\\
        { - 1}&{6 - {k^2}{h^2}}&{ - 1}\\
        0&{ - 1}&0
    \end{array}} \right]_{x_3}
    \left[ {\begin{array}{*{20}{c}}
        0&0 &0\\
        0&-1&0\\
        0&0 &0
    \end{array}} \right]_{x_3+h}
    \right].
    \end{equation}
    
    In case the Sommerfeld radiation condition \eqref{eq:SommerfeldBC} is used, the discrete schemes for the boundary points need to be defined. We can introduce ghost points located outside the boundary points. For instance, suppose $u_{0,i_2,i_3}$ is a ghost point on the left of $u_{1,i_2,i_3}$, the normal derivative can be approximated  by
    \begin{equation}
        \frac{\partial u}{\partial \vec{n}} -\text{i} k u =\frac{u_{0,i_2,i_3}-u_{2,i_2,i_3}}{2h}-\text{i} k u_{1,i_2,i_3}= 0. \nonumber
    \end{equation}
    We can rewrite it as
    \begin{equation}\label{eq: elimination}
        u_{0,i_2,i_3} = u_{2,i_2,i_3} + 2h\text{i} k u_{1,i_2,i_3}.
    \end{equation}
    Then one can eliminate the ghost point in the stencil for the boundary points.
    
    In addition, for the discretization of the Dirac function $\delta \left( {x_1-x_0,x_2-y_0,x_3-z_0} \right)$, we can set the right-hand side (RHS) as
    \begin{equation}\label{eq: discrete dirac func}
        b_h \left( {i_1,i_2,i_3} \right) 
        =\left\{ {\begin{array}{*{20}{c}}
        { \frac{1}{h^3} \;\;\;{x_1}_{i_1} = x_0,\;{x_2}_{i_2} = y_0,\;{x_3}_{i_3} = z_0}\\
        {0\;\;\;\;\;\;\;{x_1}_{i_1} \ne x_0,\;{x_2}_{i_2} \ne y_0,\;{x_3}_{i_3} \ne z_0}
        \end{array}} \right. . \nonumber
    \end{equation}
    
    The finite-difference discretization of the partial equation on geometric grids results in a system of linear equations 
    \begin{equation}\label{eq:linearsys}
    A_h \mathbf{u}_h=\mathbf{b}_h. \nonumber
    \end{equation}
    With Sommerfeld boundary  conditions, the resulting matrix is sparse, symmetric, complex-valued, indefinite, and non-Hermitian for a sufficiently large wavenumber $k$.
    
    Note that $kh$ is an important parameter that can indicate how many grid points per wavelength are used.
    The grid width  $h$ can be determined by the rule of thumb of including at least $N_{pw}$(e.g., 10 or 30) grid points per wavelength. One has the following relationships 
    \begin{equation}
    kh=\frac{2 \pi h}{\lambda} = \frac{2 \pi }{N_{pw}}. \nonumber
    \end{equation}
    For example, if at least 10 grid points per wavelength are required, one has to satisfy the condition $kh \leq 0.625$.

\section{Preconditioned Krylov method} \label{sec:methods}
    For a large, sparse system matrix $A$, Krylov subspace methods are popular. This section will specify every component of the preconditioned Krylov subspace methods that we use to solve the linear system.
    
    \subsection{Krylov subspace methods}
    The Krylov subspace methods are established on a collection of iterants in the subspace 
    \begin{equation}
    {K^k}\left( {A;{\mathbf{r}^0}} \right): = {\text{span}}\left\{ {{\mathbf{r}^0},A{\mathbf{r}^0}, \ldots ,{A^{k - 1}}{\mathbf{r}^0}} \right\}, \nonumber
    \end{equation}
    where $K^k$ is at most a $k$-dimensional Krylov space with respect to matrix $A$ and initial residual $\mathbf{r}^0$. For a basic iterative method, after $k$ iterations, $\mathbf{u}^k$ will be an element of $\mathbf{u}^0 + K^k (A,\mathbf{r}^0)$. Some representative Krylov methods, like Conjugate Gradient (CG) \cite{hestenes1952methods}, CGNR \cite{saad2003iterative}, MINRES \cite{paige1975solution}, BICG \cite{fletcher1976conjugate}, Bi-CGSTAB \cite{van1992bi}, GMRES \cite{saad1986gmres}, are developed so far. Among these, the CG method is a basic one. The error is minimized in the A-norm, and it only needs three vectors in memory during iterations. However, this algorithm is mainly designed for the symmetric and positive definite system matrix. In contrast, Bi-CGSTAB and GMRES can be used for consistent problems that are indefinite and non-symmetric.
    They are suitable choices for the Helmholtz equation. Compared with full GMRES, Bi-CGSTAB has short recurrences and better parallel properties.

    Also, the IDR(s),  which is developed by Sonneveld, Van Gijzen \cite{sonneveld2009idr}, is an efficient alternative to Bi-CGSTAB for Helmholtz problems \cite{knibbe2011gpu}. In IDR(s), $s$ pre-defined vectors are used to enhance the convergence. \cite{sonneveld2009idr} showed that IDR(1) has similar computational complexity and memory requirements as Bi-CGSTAB. With higher values of $s$, IDR(s) show performance close to GMRES but with more storage requirements. For example, we need to store 17 vectors for IDR(4), while Bi-CGSTAB needs to store 7 vectors. In general, IDR(4) is sufficient for most of the problems. Thus, we employ IDR(4) as a representation of IDR(s) in our parallel framework. For a detailed analysis of the influence of $s$ on the algorithm's convergence properties, we refer the readers to the work of  Sonneveld and Van Gijzen et al. \cite{sonneveld2009idr,Collignon2011Minimizing}
    
    We first focus on implementing the matrix-free parallelization of GMRES, then it can be directly generalized to Bi-CGSTAB and IDR(s).
    
    \subsection{Preconditioning} \label{sec:lrpre}
    The coefficient matrix of a non-Hermitian linear system is preferred to have a spectrum located in a bounded region that excludes the origin in the complex plane, which results in fast convergence for iterative methods. We can incorporate a preconditioner to enhance the convergence of Krylov subspace methods. That is to pre-multiply the linear system with a preconditioning matrix $M_h^{-1}$. It can be implemented by left preconditioning
    \begin{equation}
    M_h^{-1}A_h \mathbf{u}_h=M_h^{-1} \mathbf{b}_h, \nonumber
    \end{equation}
    or right preconditioning,
    \begin{equation}
    A_h M_h^{-1} \tilde{\mathbf{u}}_h=\mathbf{b}_h, \nonumber
    \end{equation}
    where $\mathbf{u}_h = M_h \tilde{\mathbf{u}}_h$. It can be proved that there is no essential difference between both preconditioning methods with respect to convergence behavior. But it is worth noting that the residual vectors computed by left preconditioning and right preconditioning correspond to the preconditioned and actual residuals, respectively. Thus, it may cause some difficulties for left preconditioning if a stopping criterion needs to be based on the actual residual instead of the preconditioned one.
    
    \subsection{Complex Shifted Laplace Preconditioner} \label{sec:cslp}
    We focus on CSLP motivated by its nice performance and easy setup. 
    The CSLP is defined by
    \begin{equation}\label{eq:preconditioner}
    M_h =-\Delta_h-\left(\beta_1-\beta_2 \text{i} \right) k^2 I_h,
    \end{equation}
    where $\beta_2$ is so-called the complex shift. It is proved that a small magnitude of the complex shift is necessary for convergence if the wavenumber $k$ is large \cite{gander2015applying}. In practice, it is also expensive to invert the preconditioner. Standard multigrid methods \cite{maclachlan2008algebraic,erlangga2006novel} are usually employed to invert a shifted Laplacian preconditioner. However, as the complex shift gets small, one cannot ensure that the standard multigrid methods will be effective \cite{cocquet2017large}. Besides, the computational burden to approximate the preconditioner will significantly increase. Thus, it is necessary to choose a proper complex shift. In the numerical experiments of this paper, unless noted otherwise, $\beta_1=1, \beta_2=-0.5$ will be used \cite{Gijzen2007} and the inverse of the preconditioner will be approximated by one geometric multigrid V-cycle.

    \subsection{3D multigrid for the preconditioner solve}
    One needs to compute the inverse of preconditioner $M$ in the preconditioned Krylov-based algorithms. A direct solver is a good choice if the preconditioner is simple, such as the block-Jacobi preconditioner. However, it is usually too costly to directly invert a preconditioner like CSLP. One idea is to approximately solve the preconditioner by using the multigrid method \cite{erlangga2006novel}. 
    
    A 3D multigrid method involves several components that need a careful design to achieve excellent convergence. The first ingredient of multigrid methods is the smoothing property of basic iterative methods. For example, the Gauss-Seidel and SOR($\omega$) methods can serve as efficient smoothers. As coarser grids are involved,  going back and forth between different hierarchies of grids is required. Specifically, the inter-grid transfer operations include restriction and prolongation operators. Besides, the construction of coarse grid operators also needs careful attention.
    
    \subsubsection{Smoothers}
    Classical iteration methods such as the Gauss-Seidel or damped Jacobi iterations can be used as smoothers. The impact of varying the relaxation parameter for the damped Jacobi smoother and using different smoothers on the convergence properties have already been extensively investigated by \cite{erlangga2005robust} and \cite{Hocking2021Optimal}. It has been shown that, based on an actual situation, the use of smoothers and the choice of relaxation parameters can be flexible within a certain range. In this framework, we will mainly use the damped Jacobi smoother as it is easy to parallelize and has been shown effective in the before-mentioned papers. We will fix the relaxation parameter at $\omega = 0.8$, as this value was found to work well for our test problems.

    \subsubsection{Restriction}
\begin{figure}[htbp]
        \centering
        \includegraphics[width=0.4\columnwidth]{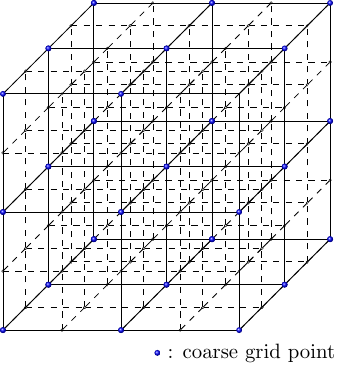}
        \caption{Vertex-centered 3D standard coarsening.}
        \label{fig: fine grid with coarse grid}
    \end{figure}

    Without loss of generality, the inter-grid operations will be presented in a simple case of two grids, i.e., a fine, and a coarse grid. Two sets of uniform grids with size $h$ and $H=2h$ are used to discretize a regular computational domain $\Omega = (0,1) \times (0,1) \times (0,1)$. Figure \ref{fig: fine grid with coarse grid} shows part of a 3D fine grid with a coarse grid obtained by standard coarsening.

    The 3D full-weighting restriction operator has a stencil given by
    \begin{equation}\label{eq:weight operator}
    I_h^H = \frac{1}{64}\left[ 
    \left[{\begin{array}{*{20}{c}}
        1&2&1\\
        2&4&2\\
        1&2&1
        \end{array}} \right]_h^H
    \left[{\begin{array}{*{20}{c}}
        2&4&2\\
        4&8&4\\
        2&4&2
        \end{array}} \right]_h^H
    \left[{\begin{array}{*{20}{c}}
        1&2&1\\
        2&4&2\\
        1&2&1
        \end{array}} \right]_h^H
    \right]. \nonumber
    \end{equation}
    Let $b^H = I_h^H r^h$, where $b^H$ is the right-hand side of the coarse grid, and $r^h$ is the residual of the fine grid. The restriction can be implemented in a matrix-free way as follows.
    \begin{equation}\label{eq: matrix-free restriction}
        \begin{aligned}
            b_{i_1,i_2,i_3}^H &= I_h^H r_{2i_1,2i_2,2i_3}^h = \frac{1}{64} \left( 8r_{2i_1,2i_2,2i_3}^h \right.\\
            & + 4r_{2i_1-1,2i_2,2i_3}^h + 4r_{2i_1+1,2i_2,2i_3}^h + 4r_{2i_1,2i_2-1,2i_3}^h \\
            & + 4r_{2i_1,2i_2+1,2i_3}^h + 4r_{2i_1,2i_2,2i_3-1}^h + 4r_{2i_1,2i_2,2i_3+1}^h \\
            & + 2r_{2i_1-1,2i_2-1,2i_3}^h + 2r_{2i_1+1,2i_2-1,2i_3}^h + 2r_{2i_1-1,2i_2+1,2i_3}^h + 2r_{2i_1+1,2i_2+1,2i_3}^h \\
            & + 2r_{2i_1-1,2i_2,2i_3-1}^h + 2r_{2i_1+1,2i_2,2i_3-1}^h + 2r_{2i_1-1,2i_2,2i_3+1}^h + 2r_{2i_1+1,2i_2,2i_3+1}^h \\
            & + 2r_{2i_1,2i_2-1,2i_3-1}^h + 2r_{2i_1,2i_2+1,2i_3-1}^h + 2r_{2i_1,2i_2-1,2i_3+1}^h + 2r_{2i_1,2i_2+1,2i_3+1}^h \\
            & + r_{2i_1-1,2i_2-1,2i_3-1}^h + r_{2i_1+1,2i_2-1,2i_3-1}^h + r_{2i_1-1,2i_2+1,2i_3-1}^h + r_{2i_1+1,2i_2+1,2i_3-1}^h \\
            & + \left. r_{2i_1-1,2i_2-1,2i_3+1}^h + r_{2i_1+1,2i_2-1,2i_3+1}^h + r_{2i_1-1,2i_2+1,2i_3+1}^h + r_{2i_1+1,2i_2+1,2i_3+1}^h \right),
        \end{aligned}\nonumber
    \end{equation}
    where $(i_1,i_2,i_3) \in \Omega^H$.

    \subsubsection{Interpolation}
    The interpolation operator $I_H^h$ transfers grid vectors from the coarse to the fine grid. The 3D trilinear interpolation operator stencil can be written as
    \begin{equation}\label{eq:bilinear operator}
    I_H^h = \frac{1}{8}
    \left[
    \left[{\begin{array}{*{20}{c}}
            1&2&1\\
            2&4&2\\
            1&2&1
            \end{array}} \right]_H^h
    \left[{\begin{array}{*{20}{c}}
            2&4&2\\
            4&8&4\\
            2&4&2
            \end{array}} \right]_H^h
    \left[{\begin{array}{*{20}{c}}
            1&2&1\\
            2&4&2\\
            1&2&1
            \end{array}} \right]_H^h
    \right]. \nonumber
    \end{equation}
    
    \begin{figure}[htbp]
        \centering
        \includegraphics[width=0.3\columnwidth]{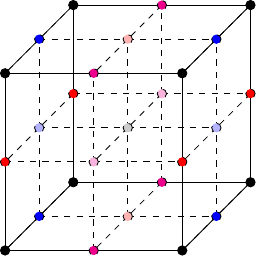}
        \caption{The allocation map of interpolation operator.}
        \label{fig: interpolation stencil}
    \end{figure}
    Let $u^h = I_H^h u^H$, where $u^H$ is the solution of the coarse grid, and $u^h$ is the correction for the fine grid. As shown in Figure \ref{fig: interpolation stencil}, the interpolation can be implemented as follows.
    
    \begin{equation}\label{eq: matrix-free prolongation}
        \begin{aligned}
            &I_H^hu_{{i_1},{i_2},{i_3}}^H = \\
            &\left\{ \begin{aligned}
            & u_{{i_1},{i_2},{i_3}}^H, &&\textcolor{black}{\bullet} \\
            & \frac{1}{2}\left( {u_{{i_1},{i_2},{i_3}}^H + u_{{i_1} + 1,{i_2},{i_3}}^H} \right), &&\textcolor{blue}{\bullet} \\
            & \frac{1}{2}\left( {u_{{i_1},{i_2},{i_3}}^H + u_{{i_1},{i_2}+1,{i_3}}^H} \right), &&\textcolor{magenta}{\bullet} \\
            & \frac{1}{2}\left( {u_{{i_1},{i_2},{i_3}}^H + u_{{i_1},{i_2},{i_3}+1}^H} \right), &&\textcolor{red}{\bullet} \\
            & \frac{1}{4}\left( {u_{{i_1},{i_2},{i_3}}^H + u_{{i_1}+1,{i_2},{i_3}}^H + u_{{i_1},{i_2}+1,{i_3}}^H + u_{{i_1}+1,{i_2}+1,{i_3}}^H} \right), &&\textcolor{red!30}{\bullet} \\
            & \frac{1}{4}\left( {u_{{i_1},{i_2},{i_3}}^H + u_{{i_1}+1,{i_2},{i_3}}^H + u_{{i_1},{i_2},{i_3}+1}^H + u_{{i_1}+1,{i_2},{i_3}+1}^H} \right), &&\textcolor{blue!30}{\bullet} \\
            & \frac{1}{4}\left( {u_{{i_1},{i_2},{i_3}}^H + u_{{i_1},{i_2}+1,{i_3}}^H + u_{{i_1},{i_2},{i_3}+1}^H + u_{{i_1},{i_2}+1,{i_3}+1}^H} \right), &&\textcolor{magenta!30}{\bullet} \\
            & \frac{1}{8}\left( u_{{i_1},{i_2},{i_3}}^H + u_{{i_1},{i_2}+1,{i_3}}^H + u_{{i_1}+1,{i_2},{i_3}}^H  \right.\\
            & + u_{{i_1}+1,{i_2}+1,{i_3}}^H + u_{{i_1},{i_2},{i_3}+1}^H + u_{{i_1},{i_2}+1,{i_3}+1}^H \\
            & + \left. u_{{i_1}+1,{i_2},{i_3}+1}^H + u_{{i_1}+1,{i_2}+1,{i_3}+1}^H \right), &&\textcolor{black!20}{\bullet} 
            \end{aligned} \right.
        \end{aligned}
    \end{equation}
    where $(i_1,i_2,i_3) \in \Omega^H$.

    \subsubsection{Coarse grid operator}
    The coarse grid matrix $A_{H}$ can be built in two ways. The first way is to obtain $A_H$ by re-discretizing on the coarse grid in the same way as the matrix $A_h$ is obtained on the fine grid, which is known as the discretized coarse grid operator (DCG). The second approach is known as Galerkin Coarsening $A_H=I_h^H A_h I_H^h$. Galerkin coarse-grid operator GCG is more general in its range of applicability, but it has an associated expense in terms of growing stencils. Considering parallel implementation, we choose DCG as our method. Boundary conditions of the preconditioner operator are set identically to the corresponding model problems.

    \subsubsection{Multigrid cycles}
	
    Based on the ingredients above, a classical two-Grid cycle is given in Algorithm \ref{two-grid cycle}, where $S_h$ denotes the smoothing procedure, and $\upsilon_1$ and $\upsilon_2$ are the numbers of pre- and post-smoothing steps, respectively.
    
    \begin{center}
        \scalebox{0.98}{
            \begin{algorithm}[H]
                \SetAlgoLined
                $\upsilon_1$ pre-smoothing sweeps: $\mathbf{u}_h^{i+1/3}=S_h^{\upsilon_1}(\mathbf{u}_h^i, A_h, \mathbf{b}_h)$\;
                Residual computation: $\mathbf{r}_h=\mathbf{b}_h-A_h \mathbf{u}_h^{i+1/3}$\;
                Restriction of the residual to $G_H$: $\mathbf{r}_H=I_h^H \mathbf{r}_h$\;
                Determination of the error on $G_H$: $\mathbf{e}_{H}=\left( A_H\right)^{-1} \mathbf{r}_H $\;
                Prolongation of the error to $G_h$: $\mathbf{e}_h= I_H^h \mathbf{e}_H$\;
                Correction of the last solution iterate: $\mathbf{u}_h^{2/3}=\mathbf{u}_h^{1/3}+\mathbf{e}_h$\;
                $\upsilon_2$ post-smoothing sweeps: $\mathbf{u}_h^{i+1}=S_h^{\upsilon_2} (\mathbf{u}_h^{i+2/3}, A_h, \mathbf{b}_h )$.
                \caption{A two-grid cycle.}
                \label{two-grid cycle}
            \end{algorithm}
        }
    \end{center}
    
    Two-grid methods are rarely practical because the coarse-grid problem may still be too large to be solved with a direct method. The idea to apply the two-grid idea to $A_H$ recursively and to continue until the coarse linear system can be solved with negligible computational costs gives rise to a genuine multigrid method.
    
    Different \textit{cycle} types can be distinguished based on the sequence in which the grids are traversed within a single multigrid iteration. The V-cycle and W-cycle are achieved by calling the two-grid method once or twice, respectively, on each coarse grid. The V-cycle method involves traversing from the finest grid to the coarsest grid, performing a single pre-smoothing step at each level, and then moving back up to the finest grid, performing a single post-smoothing step at each level. The F-cycle lies between the V-cycle and the W-cycle. It begins with the restriction to the coarsest grid. During the prolongation process, after reaching each level for the first time, an additional V-cycle to the coarsest grid is performed. In the numerical experiments, both V-cycle and F-cycle methods will be involved.
    
    For solving the coarse grid problems (i.e., step 4 in Algorithm \ref{two-grid cycle}), a direct solver is not easy to be parallelized in a matrix-free way. Thus, we solve the final coarse-grid problem iteratively by GMRES. 

    \subsection{Matrix-free method}\label{Matrix-free Matvec}
    We can implement the Krylov subspace methods in a matrix-free way instead of constructing the coefficient matrices explicitly. The matrix-vector multiplication can be replaced by stencil computations. Likewise, the preconditioning matrix $M$ and its stencil can be obtained analytically by its definition \eqref{eq:preconditioner}. One does not need to construct $M$ explicitly since the result of $Mx$ can be calculated by its corresponding stencil. 
    
    Considering any grid point $(i_1, i_2, i_3)$, define $ap$, $aw$, $ae$, $as$, $an$, $ad$ and $au$ as the multipliers of $u(i_1, i_2, i_3)$, $u(i_1-1, i_2, i_3)$, $u(i_1+1, i_2, i_3)$, $u(i_1, i_2-1, i_3)$, $u(i_1, i_2+1, i_3)$, $u(i_1, i_2, i_3-1)$ and $u(i_1, i_2, i_3+1)$, respectively. When physical boundary conditions are encountered, it only needs to eliminate the ghost grid points as given in equation \eqref{eq: elimination}.
    
    For the Helmholtz operator \eqref{eq:Helmstencil}, we have
    \begin{equation}
        ap = \frac{6-k^2 h^2}{h^2} \quad  aw = ae = as = an = au = ad= -\frac{1}{h^2}. \nonumber
    \end{equation} 
    
    As we invert the preconditioner $M$ by the multigrid method, $y = M_h x$ and $y = M_H x$ need to be performed in the smoother and coarsest grid solver. For these CSLP operators, according to equation \eqref{eq:stencilLaplace} and \eqref{eq:preconditioner}, we will have
    \begin{equation}
        ap = \frac{6-\left(\beta_1 - \beta_2\text{i}  \right) k^2 h^2}{h^2} \quad  aw = ae = as = an = au = ad = -\frac{1}{h^2}. \nonumber
    \end{equation}
    
    The results of $v=A_h u$, $y = M_h x$ and $y = M_H x$ can be computed in a matrix-free way by Algorithm \ref{Matrix free Mat-vec}.
    
    Besides, the matrix-free restriction and prolongation operators in the multigrid method can be implemented according to equations \eqref{eq: matrix-free restriction} and \eqref{eq: matrix-free prolongation}, respectively.
    
    \begin{center}
        \scalebox{0.98}{
            \begin{algorithm}[H]
                \SetAlgoLined
                \KwIn{3D Array $u(1:n_1,1:n_2,1:n_3)$ }
                Initiate $ap$, $aw$, $ae$, $as$, $an$, $ad$, and $au$  \;
                Initiate 3D Array $v(1:n_1,1:n_2,1:n_3)$\;
                Exchange the interface boundaries data of $u$ \;
                \For{$i_3:=1, 2, ..., n_3$}{
                \For{$i_2:=1, 2, ..., n_2$}{
                    \For{$i_1:=1, 2, ..., n_1$}{
                        \If{physical boundary grid point}{
                        Reset $ap$, $aw$, $ae$, $as$, $an$, $ad$, and $au$\;}
                        $v(i_1,i_2,i_3)=ap*u(i_1,i_2,i_3)+ae*u(i_1+1,i_2,i_3)+aw*u(i_1-1,i_2,i_3)+an*u(i_1,i_2+1,i_3)+as*u(i_1,i_2-1,i_3)+au*u(i_1,i_2,i_3+1)+ad*u(i_1,i_2,i_3-1)$\;
                    }
                }}
                Return $v$.
                \caption{Matrix-free Matrix-Vector Multiplication.}
                \label{Matrix free Mat-vec}
            \end{algorithm}
        }
    \end{center}

    \section{Parallel implementation} \label{sec:implem}
    A parallel Fortran 90 code is developed to solve the 3D heterogeneous Helmholtz problems. The MPI standard is employed for data communications among the processes. Therefore, the design of an MPI topology will be the basis. Further, the domain partition and the data structure within the processes will determine the implementation of the matrix-vector multiplication and dot product. Finally, the parallelization of the Krylov subspace methods and the multigrid cycle are carried out.
    
    \subsection{Parallel setup}
\begin{figure}[htbp]
        \centering
        \includegraphics[width=0.7\columnwidth,trim=2 60 2 2,clip]{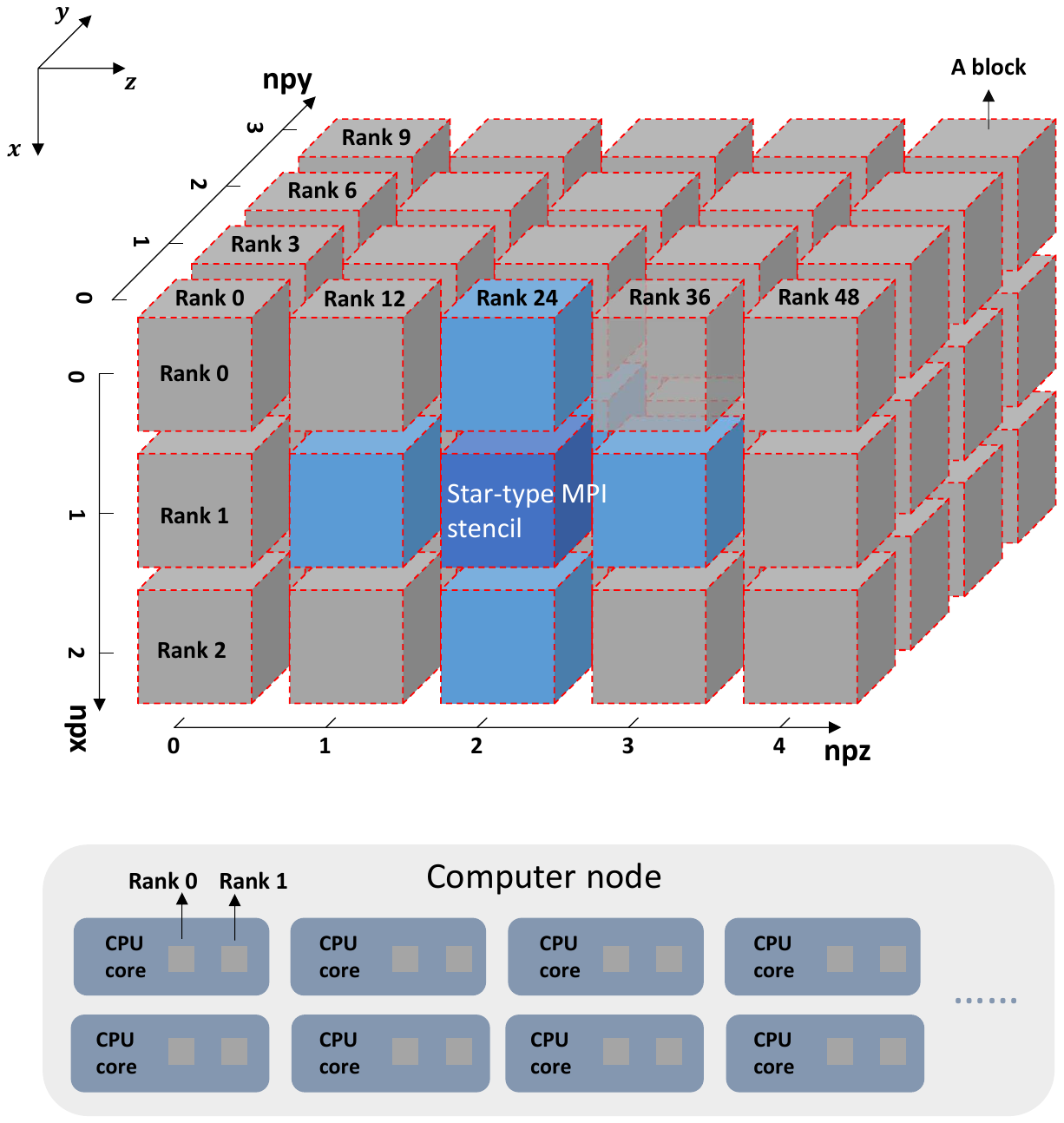}
        \caption{The schematic diagram of the MPI topology for npx0 $\times$ npy0 $\times$ npz0 = $3 \times 4 \times 5$.}
        \label{fig: MPI topology}
    \end{figure}

    For the MPI setup, the first step is to determine the total number of processes (denotes as $np$) and the number of processes in each direction. In 3D cases, the number of processes in the $x$-, $y$- and $z$- direction is denoted as $npx0$, $npy0$and $npz0$ respectively. As shown in Figure \ref{fig: MPI topology}, the entire computational domain is partitioned into $3 \times 4 \times 5$ blocks. Each block is seen as an MPI process, each of which is assigned to a CPU core. Due to the star-type computation stencil, the communication between processes exists only between the adjacent blocks in each coordinate direction. After acquiring $np$ processes and creating a parallel computing environment using MPI, each process will own its MPI rank (0 to $np-1$). We assign each subdomain to each process according to the $x$-line lexicographic order. Indicating to  different processes whether their corresponding subdomains contain physical boundaries or interface boundaries, we define $npx$, $npy$, and $npz$ to describe the position in the $x$-, $y$- and $z$- directions for every process. For example, $npx$ is in the range of $[0, npx0-1]$. When $npx$ is equal to 0, the process needs to deal with the left physical boundary. When $npx$ equals $npx0-1$, the process needs to handle the right physical boundary. 
    
\begin{figure}[htbp]
        \centering
        \includegraphics[width=0.8\columnwidth]{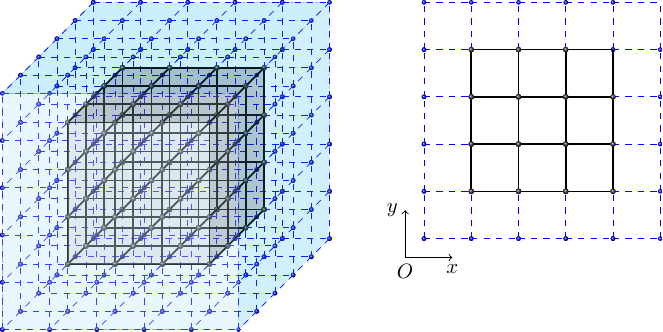}
        \caption{The schematic diagram of the local grid structure.}
        \label{fig: data structure}
    \end{figure}

    In domain partitioning, we choose to partition the domain between two grid points (i.e., along the red dotted line shown in Figure \ref{fig: global_multigrid}). Therefore, the boundary points of adjacent subdomains are adjacent grid points in the global grids. Note that the Helmholtz operator and the CSLP have a similar stencil that only needs the data from the adjacent grid points. Thus, for each block, we introduce one layer of overlapping grid points outward at each interface boundary to represent the adjacent grid points (i.e., the blue grid in Figure \ref{fig: data structure}). 
    
    In the program, the grid unknowns are stored as an array based on the grid ordering $(i_1, i_2, i_3)$ instead of a column vector based on $x$-line lexicographic ordering. We store the number of grid points in each dimension within each subdomain as $nx$, $ny$, and $nz$ respectively. To store and use the data from adjacent processors, the local arrays are extended based on the subdomain grid structure as shown in Figure \ref{fig: data structure}. Thus, the indices range of the array $\mathbf{u}$ becomes  $(1-lap: nx+lap, 1-lap: ny+lap)$. For second-order finite-difference discretization, the number of overlapping grid points $lap$ is $1$. Within a certain subdomain, the operations and array updates are limited to the range $(1:nx,1:ny,1:nz)$. The data $u(i_1,i_2,i_3)$ for $i_1 = 1$ or $nx$, $i_2=1$ or $ny$ and $i_3 = 1$ or $nz$ are sent to adjacent processors. The data received from adjacent processors are stored in the corresponding extended grid points, which are called during the operations concerning interface grid points.

    \subsection{Parallel multigrid method based on global grid}
\begin{figure}[htbp]
        \centering
        \includegraphics[width=0.7\columnwidth,trim=2 60 2 2, clip]{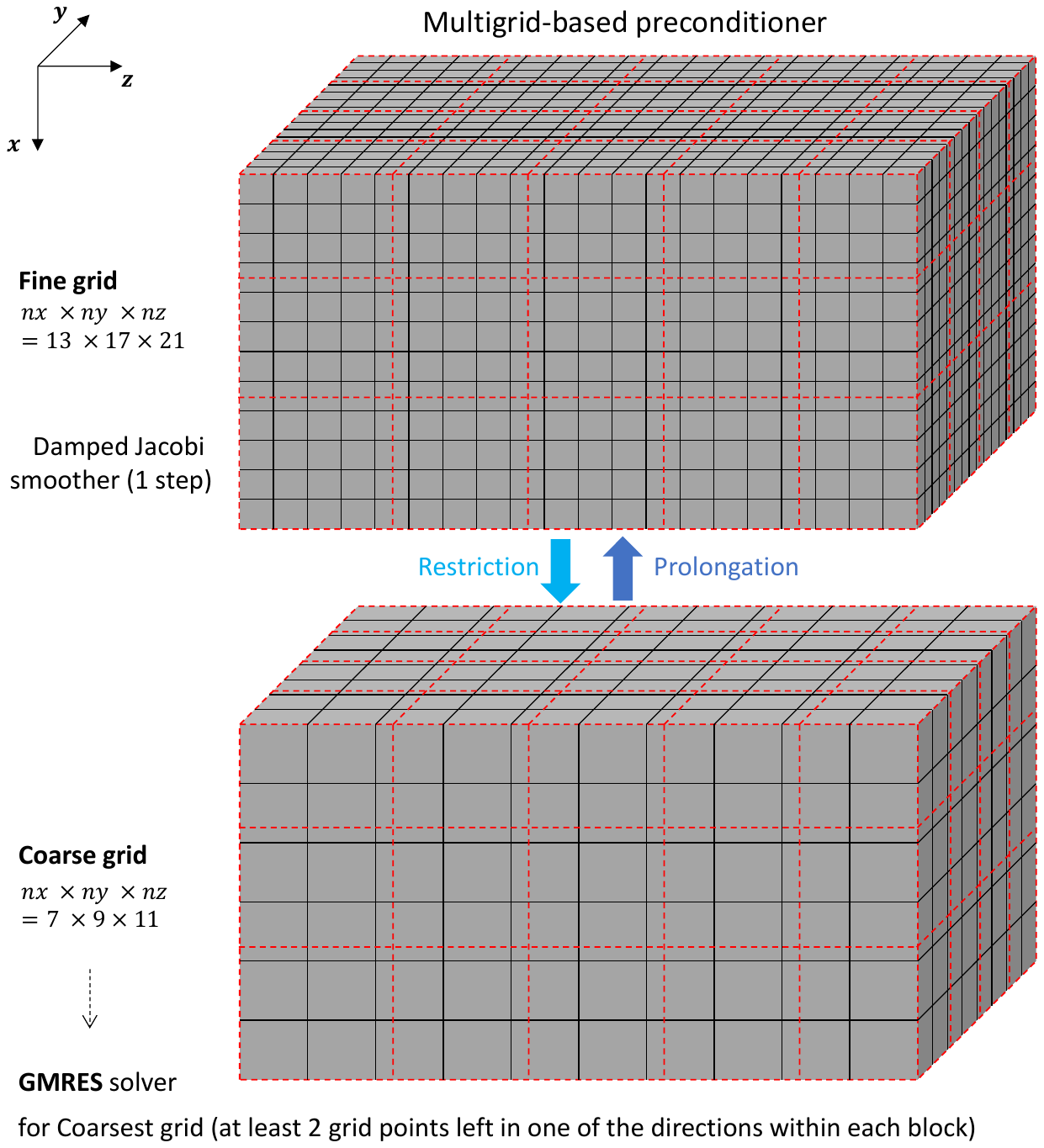}
        \caption{The schematic diagram of parallel multigrid based on global-grid coarsening. Grid size is chosen at random for demonstration purposes.}
        \label{fig: global_multigrid}
    \end{figure}

    We will consider the parallel implementation of the multigrid iteration based on the original global grid, as shown in Figure \ref{fig: global_multigrid}, where an arbitrary grid size is chosen for demonstration purposes.
    
    According to the relationship between the fine grid and the coarse grid, the parameters of the coarse grid are determined by the grid parameters of the fine one. For example, point $(i_{1c},i_{2c},i_{3c})$ in the coarse grid corresponds to point $(2i_{1c} -1, 2i_{2c} -1, 2i_{3c} -1)$ in the fine grid. The restriction, as well as interpolation of the grid variables, can be implemented according to the stencils based on the index correspondence between the coarse and fine grid. The grid-based matrix-free multigrid preconditioner starts at the finest grid, and recursively performs a two-grid cycle until the bottom level is reached, where a certain number of pre-defined grid points are left in one of the directions (denote as the stopping criterion for coarsening, \textit{i.e.} $ nx_{\text{coarsest}} \times ny_{\text{coarsest}} \times nz_{\text{coarsest}}$). The coarsest-grid problem is solved by parallel GMRES. On the one hand, direct solvers are not easy to parallelize and an exact solution for the coarse grid problem is not necessary for our preconditioner. On the other hand, we developed this work for massively parallel computing of large-scale practical problems. Agglomerating to one process and using a direct solver is not an economical option, as the all-gather communication would be expensive and the rest of the processes are idle. Thus, our idea is to stop coarsening at a certain level, on which the communication load is still much less than the computational load and we can solve the coarse-grid problem in parallel with fairly good efficiency. In this way, we can take advantage of the parallelism available in the system and achieve better overall performance. Because it allows for more accurate solutions at coarser levels, which can then be used to inform the solutions at finer levels. This is particularly important for practical applications where we need to balance both performance and accuracy. However, the optimal selection of the coarsest grid size depends on both the number of iterations and processes. Due to our limitations in computing resources and machine time, we set the stopping criterion for coarsening as $17\times 17\times 17$ to ensure that each process has at least 1-2 coarsest grid points. It should be noted that if the computational domain is not a cube, this setting means that coarsening will stop when the number of grid points in a certain direction is less than $17$. This setup allows a relatively practical lower bound for reference. The related topic of optimal coarsest grid size will be further investigated in our future work. 
    
     \section{Numerical experiments} \label{sec:exper}
    The numerical experiments are primarily conducted on the Linux supercomputer DelftBlue \cite{DHPC2022}. DelftBlue runs on the Red Hat Enterprise Linux 8 operating system. Each compute node is equipped with two Intel Xeon E5-6248R processors with 24 cores at 3.0 GHz, 192 GB of RAM, a memory bandwidth of 132 GByte/s per socket, and a 100 Gbit/s InfiniBand card. In our experiments, the solver is developed in Fortran 90. On DelftBlue, the code is compiled using GNU Fortran 8.5.0 with the compiler options \verb|-O3| for optimization purposes. Open MPI library (version 4.1.1) is employed for message passing. HDF5 1.10.7 is used for massively parallel I/O.

    For the iterative methods we used in this section, the number of matrix-vector multiplications (denote as \#Matvec), which only includes matrix-vector products with the system matrix $A_h$,  will be the main amount of work. Unless mentioned, the following convergence criterion of the preconditioned GMRES algorithm is used.
    \begin{equation}\label{eq:prcond rel res}
        \frac{\left\| M_h^{-1}\mathbf{b}_h - M_h^{-1}A_h \mathbf{u}_h^{k} \right\|_2}{\left\| M_h^{-1}\mathbf{b}_h \right\|_2} \leq 10^{-6}. \nonumber
    \end{equation}
    The convergence criterion for the preconditioned Bi-CGSTAB and IDR(s) is 
    \begin{equation}\label{eq:prcond rel res idrs}
        \frac{\left\| \mathbf{b}_h - A_h \mathbf{u}_h^{k} \right\|_2}{\left\| \mathbf{b}_h \right\|_2} \leq 10^{-6}. \nonumber
    \end{equation}
    The different stopping criteria for GMRES, Bi-CGSTAB, and IDR(4) are chosen in relation to the left and right preconditioning strategies discussed in Section \ref{sec:lrpre}. It is intended to illustrate the flexibility of our framework in handling various solvers, and preconditioning techniques. 
    
    For the coarsest-grid problem solver, according to our pre-experiments, the stopping criterion should be $3$ orders of magnitude smaller than the stopping criterion for the outer iteration. This work employs full GMRES to reduce the relative residual to $10^{-11}$, ensuring an accurate approximation for $M_H^{-1}$. Since it is only performed on the coarsest grid, it will not lead to a high computational cost.
    
    The Wall-clock time for the preconditioned Krylov solver to reach the stopping criterion is denoted as $t$. The speedup $S_p$ is defined by
    \begin{equation}
        S_p=\frac{t_r}{t_p}, \nonumber
    \end{equation}
    where $t_r$ and $t_p$ are the Wall-clock times for reference and parallel computations, respectively. The parallel efficiency $E_p$ is given by 
    \begin{equation}
        E_p=\frac{S_p}{np} =\frac{t_r}{t_p \cdot np}  , \nonumber
    \end{equation}
    where $np$ is the number of processors. It should be noted that when performing computations within a single compute node, the reference time is the wall-clock time of sequential computation. When performing distributed computing across multiple compute nodes, the reference time will be the wall-clock time of computation on a single fully-utilized node.
    
    In this section, we first validate the numerical accuracy and algorithmic flexibility of our parallel framework for three test problems. To evaluate the performance of our solver, we then focus on weak scaling, which examines the solver's performance when the problem size and the number of processing elements increase proportionally. This allows us to assess our parallel solver's ability to solve large-scale heterogeneous Helmholtz problems with minimized pollution error. We next extensively investigate the strong scaling of our parallel solution method, examining the speedup and parallel efficiency as the number of processing elements increases while the problem size remains constant. Additional performance analysis for each part of the parallel framework can be found in Appendix \ref{appendix_POP}.

    \subsection{Validation}
    We validate the numerical results by comparing numerical results from both serial and parallel computations with analytical solutions, as well as by observing the wave propagation patterns of model problems with non-constant wavenumber. Additionally, we conducted a preliminary exploration of the flexibility of various Krylov-based iterative methods and multigrid cycle types in the parallel framework.

    \subsubsection{Sequential and parallel computing}
\begin{figure}[htbp]
        \centering
        \includegraphics[width=0.49\textwidth,trim=10 10 10 10,clip]{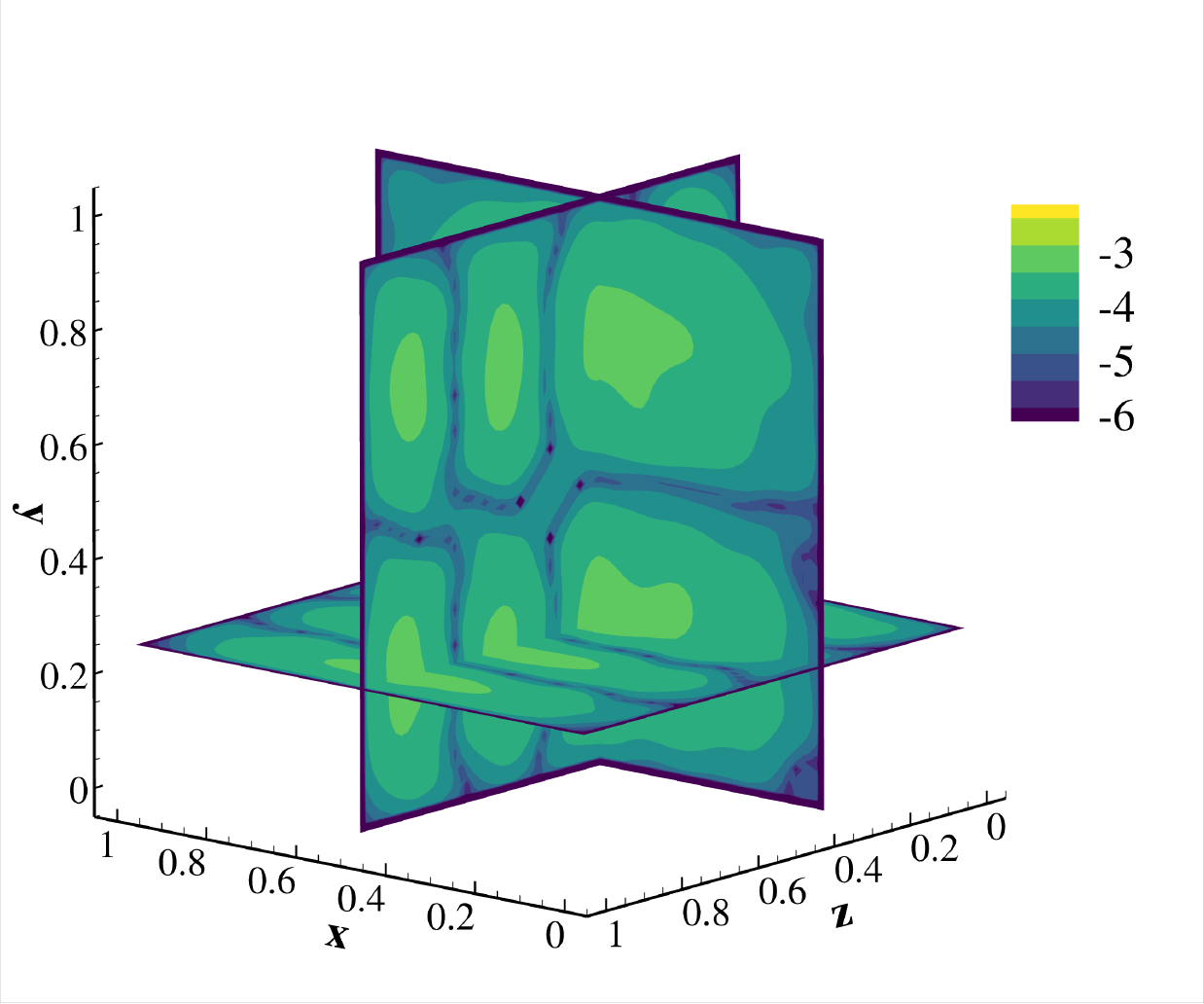}
        \hfill
        \includegraphics[width=0.49\textwidth,trim=10 10 10 10,clip]{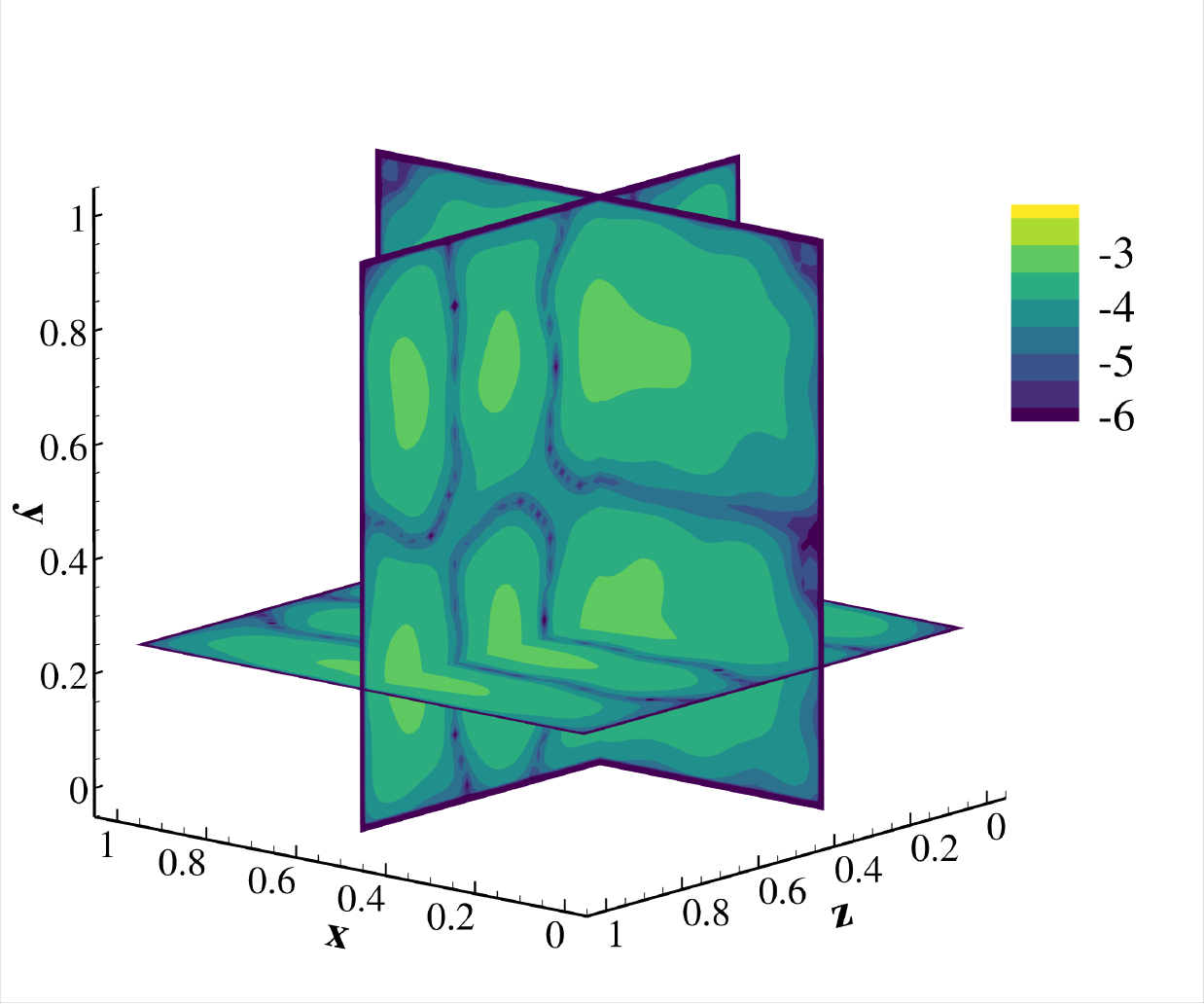}
        \caption{ The logarithmic error $\log_{10}(|u-u_{\text{analy}}|)$ of numerical results obtained by CSLP-preconditioned IDR(4) with $np=1$ (left) and $np = 3\times3\times3$ (right) for 3D Closed-off problem grid size $65 \times 65 \times 65$, $k=40$.}
        \label{fig: MP1 sol and error}
    \end{figure}

    For the accuracy validation, our parallel solver is used to solve the so-called 3D closed-off problem which has an analytical solution. According to the analytical solution given by equation \eqref{eq: MP1 analytical solution}, we can estimate the magnitude of the error between the numerical and analytical solutions. Figure \ref{fig: MP1 sol and error} shows the logarithmic error $\log_{10}(|u-u_{\text{analy}}|)$ of numerical results obtained by CSLP-preconditioned IDR(4). The left one is the result obtained by sequential computing, and the right one is the result obtained by parallel computing. Within the error tolerance, the numerical results agree well with the analytical solution, and the parallel computing setup does not introduce additional errors. It should be noted that the small but visible differences arise due to the accumulation of the rounding error. For example, the dot product $[a_1\ a_2\ a_3\ a_4]\cdot[b_1\ b_2\ b_3\ b_4]'$ calculated by $(a_1b_1+a_2b_2+a_3b_3+a_4b_4)$ and $(a_1b_1+a_3b_3)+(a_2b_2+a_4b_4)$ may yield small differences, especially for massively parallel computing. Additionally, in the aspect of the matrix, blockwise domain partitioning actually performs some elementary row and column transformations on the matrix compared to sequential computing, leading to further differences in the results. But extensive numerical experiments indicate that these differences are within an acceptable error range and will not significantly affect the convergence and accuracy.

    \subsubsection{Model problems with non-constant wavenumber}
\begin{figure}[htbp]
        \centering
        \includegraphics[width=0.45\textwidth,trim=2 2 2 2,clip]{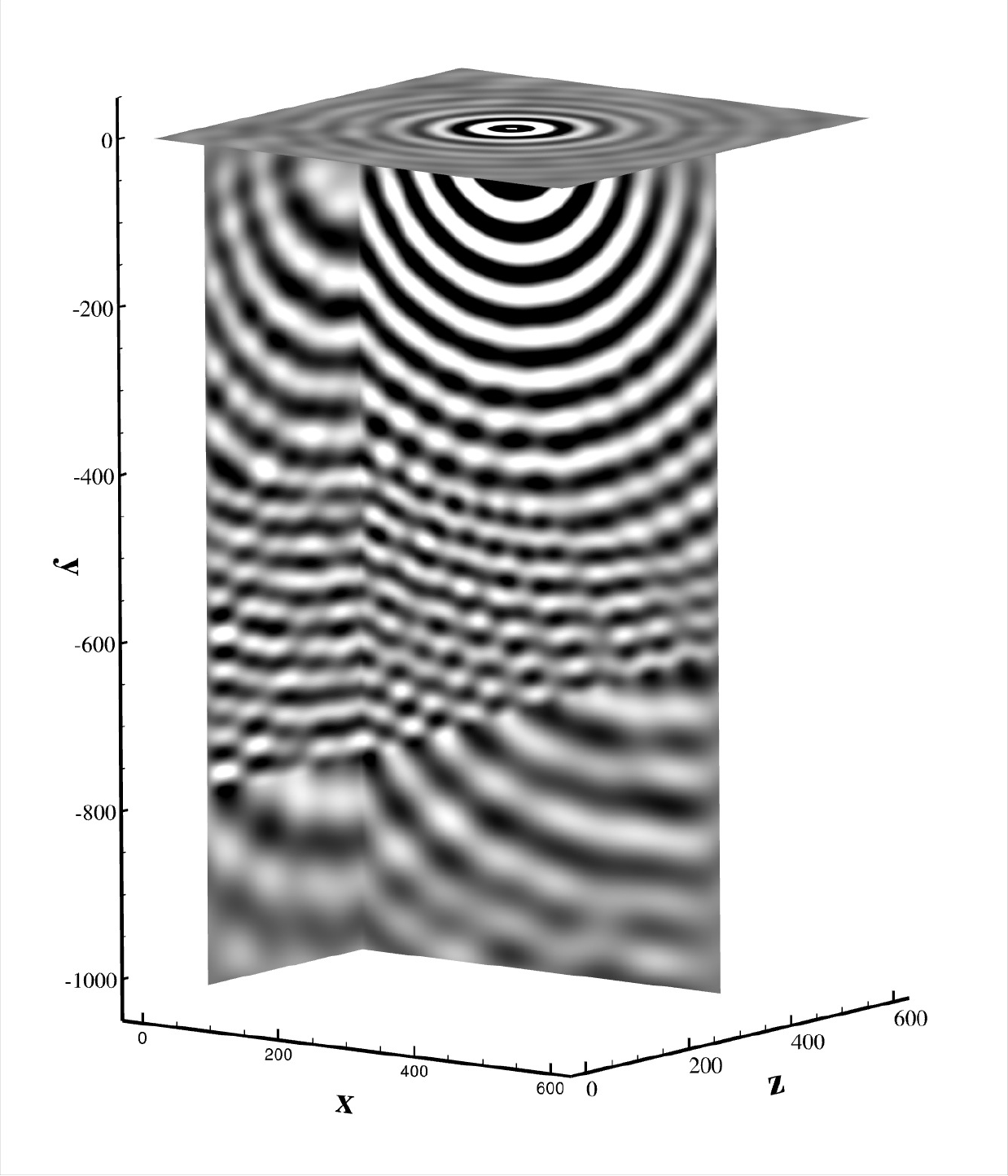}
        \hfill
        \includegraphics[width=0.45\textwidth,trim=2 2 2 2,clip]{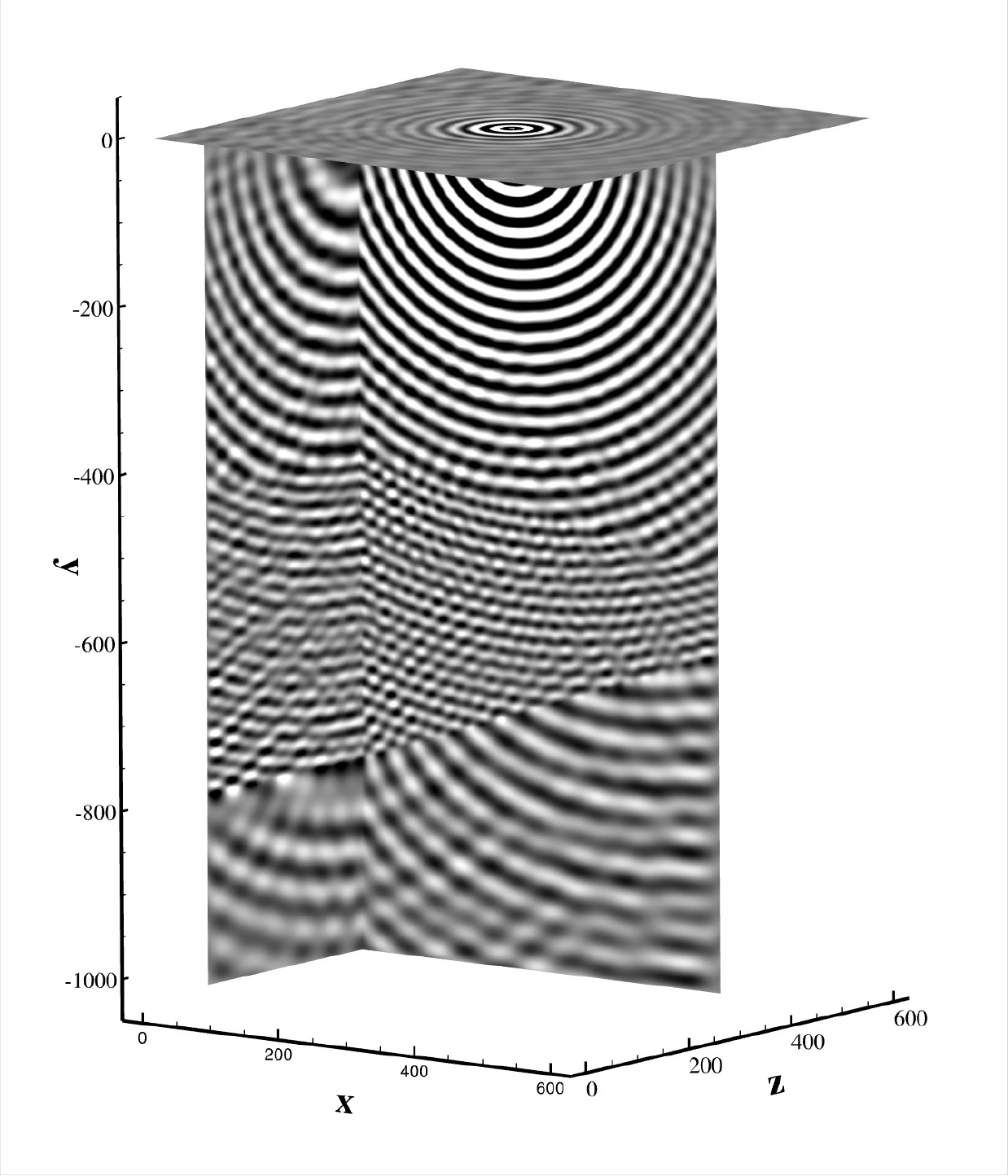}
        \caption{Real part of numerical solutions for 3D Wedge Problem with grid size $193\times 321 \times 193$ at $f=\SI{40}{\hertz}$ (left), and $385\times 641 \times 385$ at $\SI{80}{\hertz}$ (right).}
        \label{fig:MP3 sol}
    \end{figure}

    The so-called Wedge problem is a typical model with non-constant wavenumbers. 
    Figure \ref{fig:MP3 sol} shows a reasonable wave diffraction pattern of the wedge problem at $f = \SI{40}{\hertz}$ and $\SI{80}{\hertz}$, which are obtained by parallel CSLP-preconditioned IDR(4) solver. Note that the wavefront is significantly curved at the slow-to-fast interface. The wave is mainly reflected in that transition region. It illustrates that the parallel framework also works for the case of non-constant wavenumbers.
    
\begin{figure}[htbp]
        \centering
        \includegraphics[width=0.85\textwidth,trim=2 2 2 2,clip]{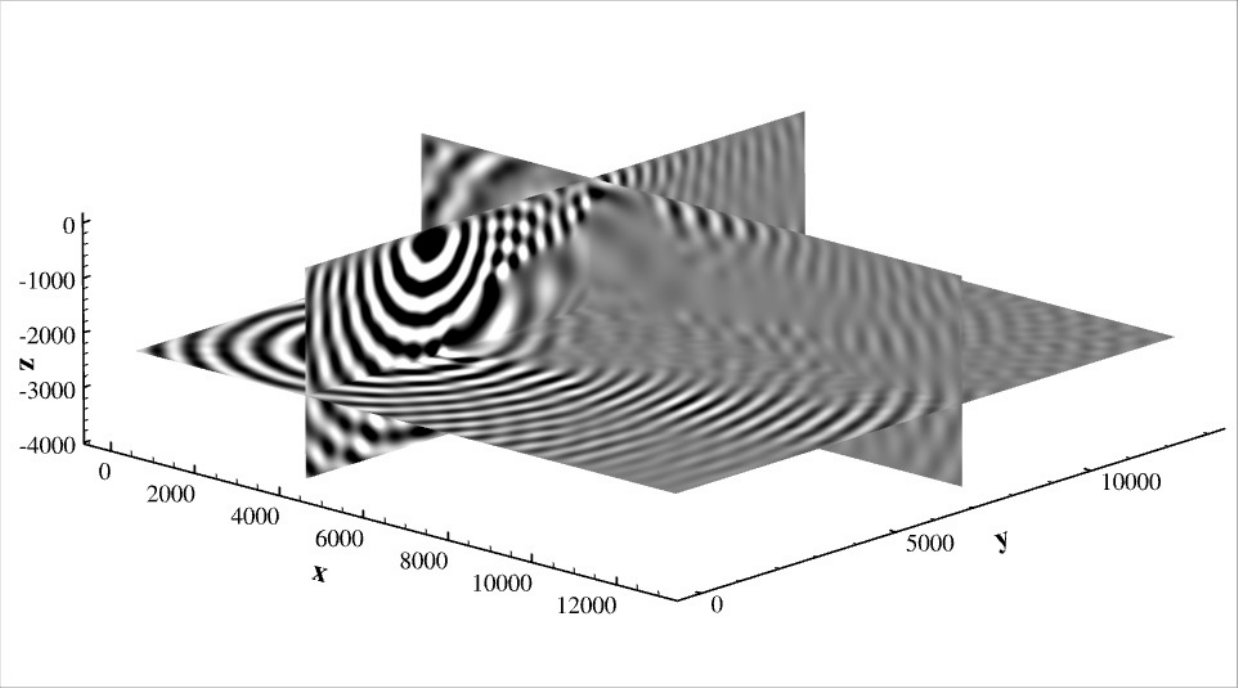}
        \caption{Real part of the solutions for 3D SEG/EAGE Salt Model with grid size $641 \times 641 \times 193$, at $f=\SI{5}{\hertz}$.}
        \label{fig:3D Salt Model Rsol}
    \end{figure}

    For practical application, we further consider the challenging 3D SEG/EAGE Salt model including complex geometries (salt domes) and a real large-size computational domain. To allow multi-level geometric coarsening, we utilize a domain with grid size $641\times641\times193$, as it has a high degree of divisibility by powers of two in each dimension. In this case, with the coarsest-grid criterion given by $17 \times 17 \times 17$,  the coarsest grid will have dimension $41 \times 41 \times 13$. A point source is located at $\left(x_1, x_2, x_3 \right) = \left( 3200, 3200, 0\right)$. Figure \ref{fig:3D Salt Model Rsol} gives a reasonable wave diffraction pattern of this model at $f=5 Hz$. It is obtained by our parallel CSLP-preconditioned IDR(4) solver. 

    \subsubsection{Krylov methods}
    Various Krylov-based iterative methods can be implemented in our parallel framework; see Table \ref{tab: matvec of diff Krlov}. For the 3D closed-off problem with Dirichlet boundary conditions, GMRES requires the least number of iterations, and IDR(4) is closer to GMRES than Bi-CGSTAB. These relative relationships can be expected in the theoretical analysis of these algorithms. One can find that the count of matrix-vector products is constant with respect to the number of processors for GMRES, while it is not for Bi-CGSTAB and IDR(4). It is because GMRES is a very robust and stable algorithm, whereas Bi-CGSTAB and IDR(s) are sensitive to the accumulation of rounding errors and the re-ordering of the linear system. Additionally, IDR(s) use random vectors to initialize, and therefore, one cannot obtain a totally constant number of iterations if the seed is not manually chosen. As Bi-CGSTAB requires a much larger number of matrix-vector multiplications, the parallel performance of GMRES and IDR(4) will be mainly considered. 

    \begin{table}[htbp]
        \centering
        \caption{The convergence behavior of different iterative methods for 3D Closed-off problem, grid size $65 \times 65 \times 65$, $k=40$.}
        \scalebox{0.9}{
          \begin{tabular}{c|c|c|c}
            \multirow{2}{*}{$np$} & \multicolumn{3}{c}{\#Matvec} \\ 
                                & GMRES  & Bi-CGSTAB  & IDR(4) \\ \hline
            1                   & 145    & 359        & 192    \\
            2$\times$1$\times$1 & 145    & 359        & 211    \\
            1$\times$2$\times$3 & 145    & 403        & 192    \\
            3$\times$3$\times$3 & 145    & 407        & 201    \\
            \end{tabular}
        }
        \label{tab: matvec of diff Krlov}
    \end{table}

    \subsubsection{Multigrid cycle types}
\begin{figure}[htbp]
        \centering
        \includegraphics[width=0.49\textwidth]{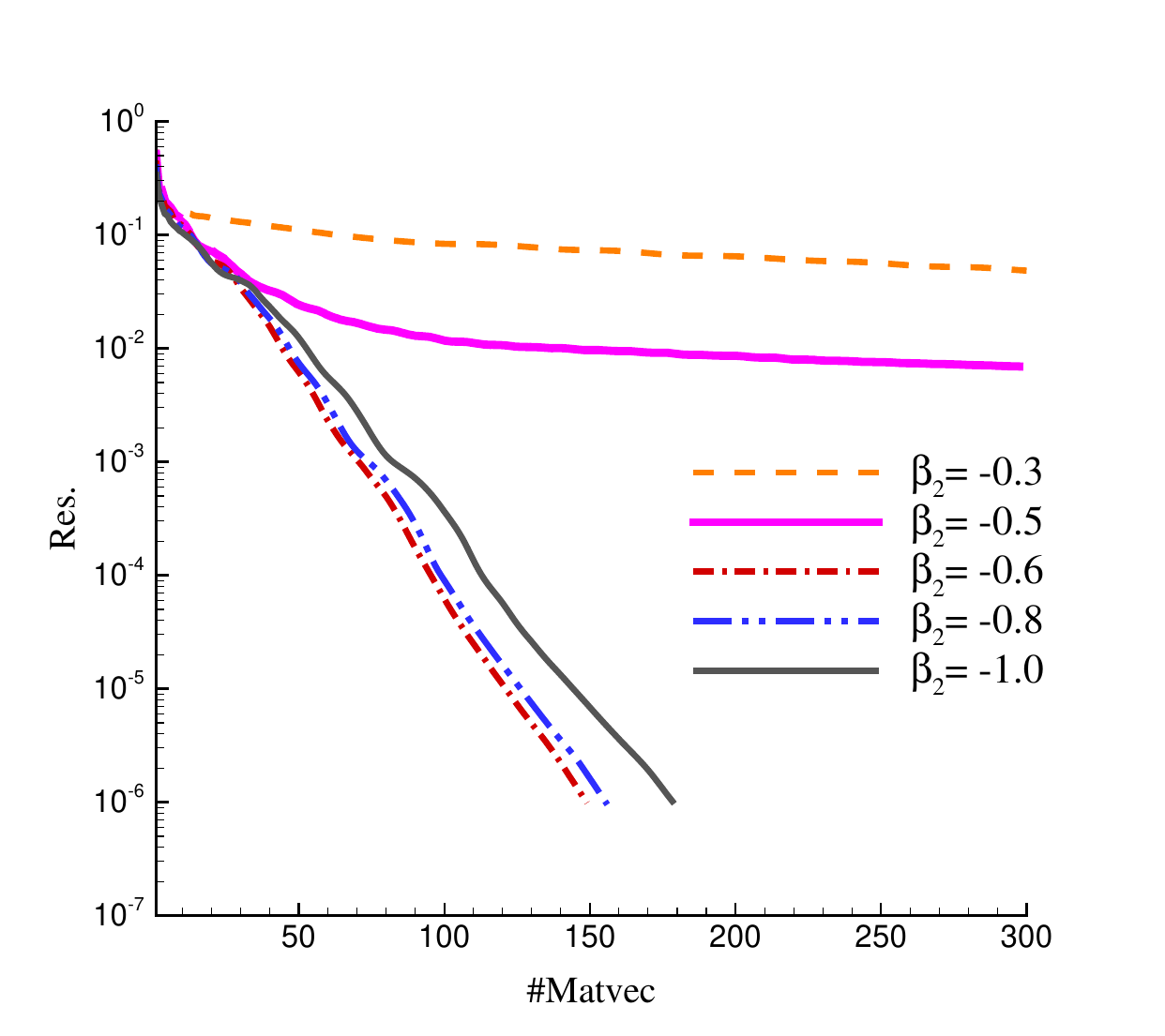}
        \hfill
        \includegraphics[width=0.49\textwidth]{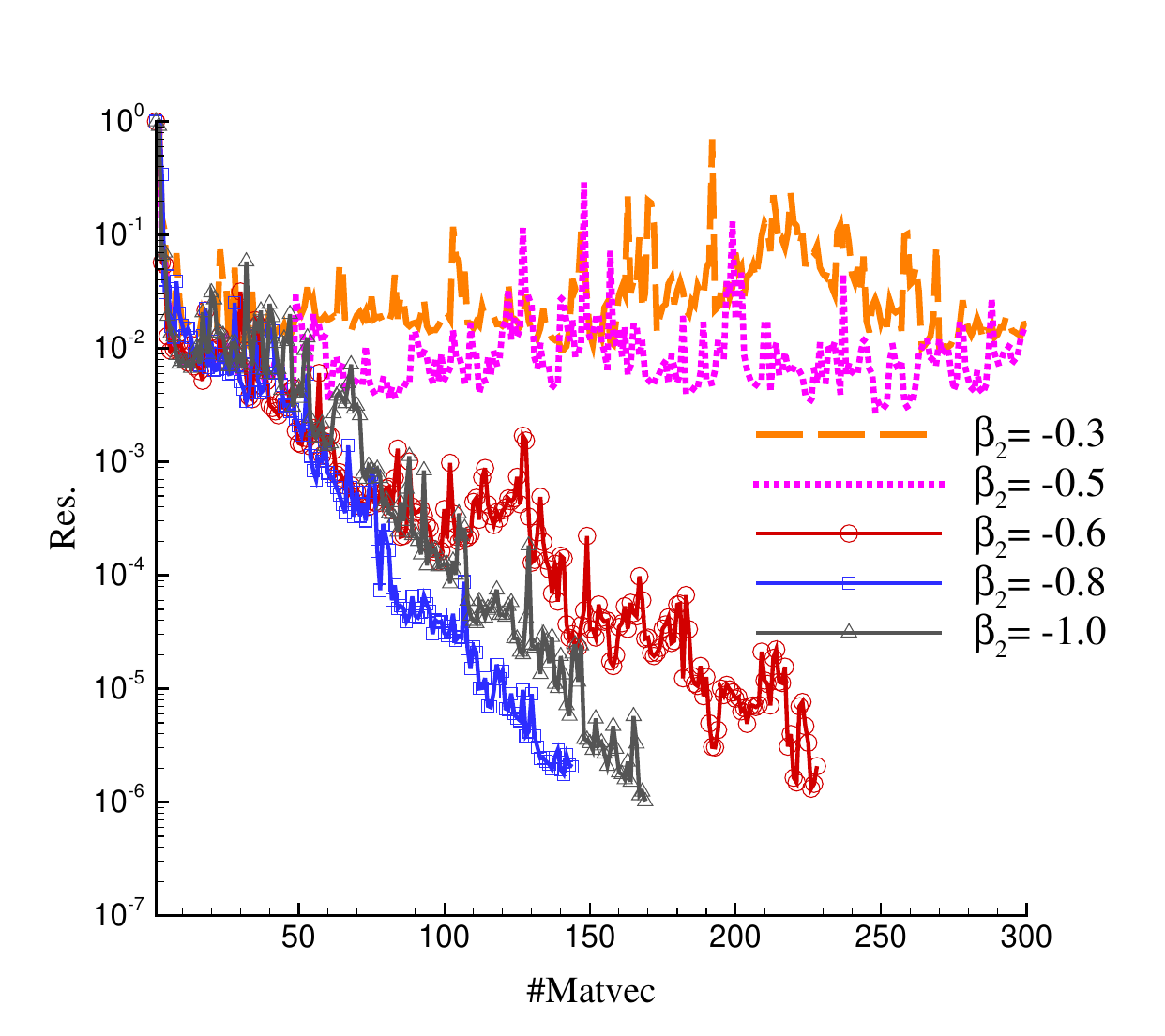}
        \caption{Convergence behavior of the parallel CSLP-preconditioned Krylov solvers, GMRES (left) and IDR(s) (right). One multigrid F-cycle is used for preconditioner solving. 3D Wedge Problem with grid size $193 \times 321 \times 193$ at $f=\SI{40}{\hertz}$ is solved. Different complex shifts of CSLP ($\beta_2$) are compared.}
        \label{fig:iters_MP3_nx193_k40_beta2s}
    \end{figure}

    Our parallel framework can solve the CSLP preconditioner not only using one multigrid V-cycle but also F-cycle. Let us take the 3D Wedge problem as an example. We found that if $\beta_2 = -0.5$, even for small frequency $f=\SI{10}{\hertz}$ with grid size $73 \times 193 \times 73$, the convergence result cannot be obtained by using the F-cycle, but it is possible by using the V-cycle. However, the Krylov solvers using one multigrid F-cycle should converge faster than those using one V-cycle theoretically. Thus, the effect of the complex shift $\beta_2$ in CSLP on convergence is studied. As shown in Figure \ref{fig:iters_MP3_nx193_k40_beta2s}, when the grid is $193 \times 321 \times 193$ and $f=\SI{40}{\hertz}$, if a larger complex shift ($-1 \le \beta_2 < -0.5$) is used, Krylov solvers can obtain convergence results by using the F-cycle. Figure \ref{fig:iters_MP3_nx193_k40_beta2s} confirms the introduction in Section \ref{sec:cslp}, namely, the choice of the complex shift parameter in the preconditioner indeed has a significant impact on the performance of multigrid methods including F-cycle as well as V-cycle. In the case of CSLP, a smaller complex shift generally leads to a better preconditioner, but the complex shift has to be large enough for standard multigrid methods to converge; see \cite{Ernst2012,gander2015applying,cocquet2017large}. For the small complex shift, moving from a V-cycle to an F-cycle may not improve the convergence. One potential explanation could be the damped Jacobi smoother used in the multigrid method in principle diverges for Helmholtz-type of problems. Therefore, using more smoothing steps, as in the F-cycle, may exacerbate the problems and lead to slower convergence compared to V-cycle. There is also evidence of this in Chapter 9 of \cite{vandanaPhD}.

\begin{figure}[htbp]
        \centering
        \includegraphics[width=0.54\textwidth]{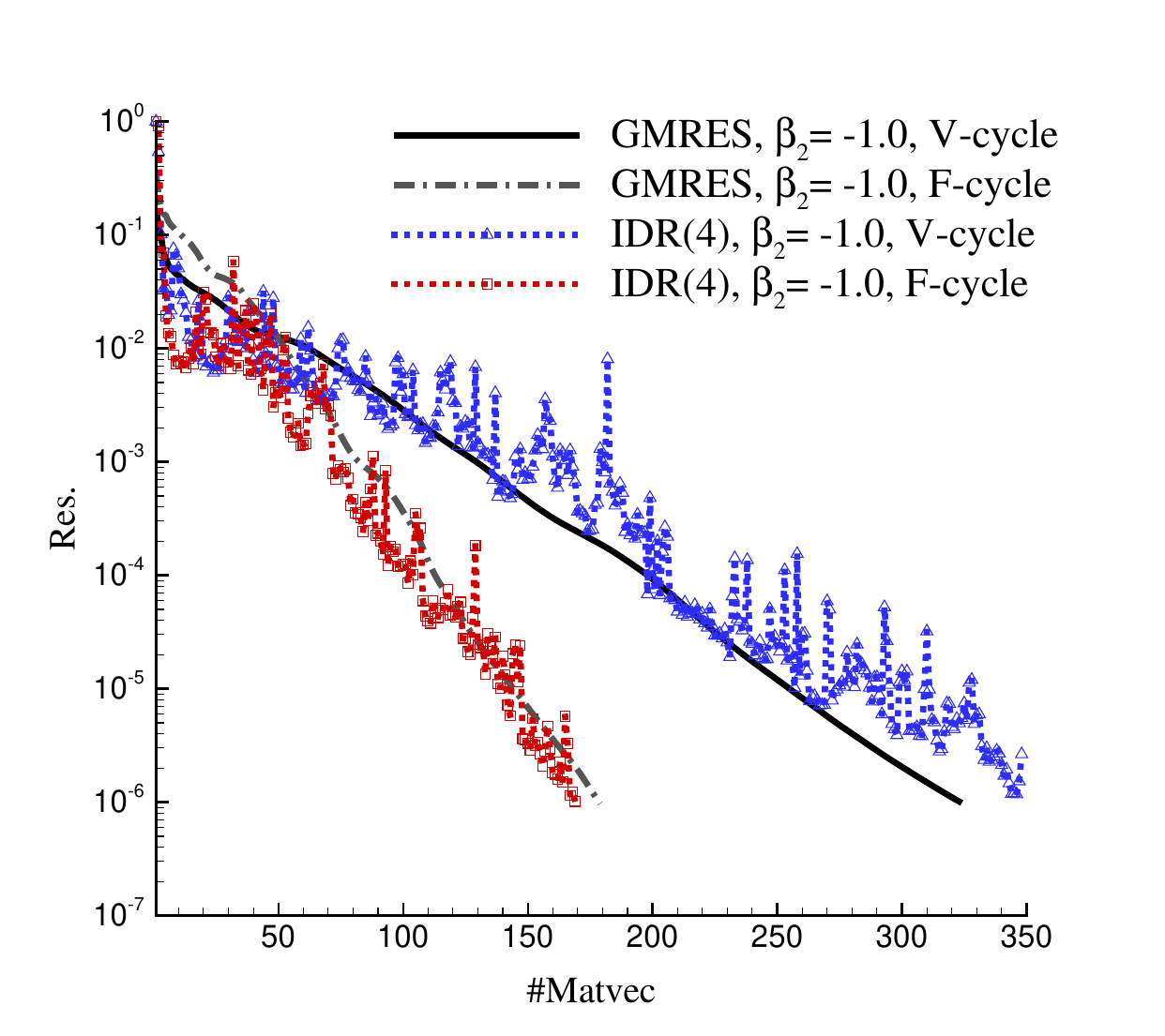}
        \caption{Convergence behavior of the parallel CSLP-preconditioned Krylov solver using one multigrid V- or F-cycle for preconditioning. The complex shift of CSLP is $\beta_2 = -1.0$. 3D Wedge Problem with grid size $193 \times 321 \times 193$ at $f=\SI{40}{\hertz}$ is solved.}
        \label{fig:iters_MP3_nx193_k40_betam1p0_idrs_gmres_F-vs-V}
    \end{figure}

    When the complex shift is $-1$, as shown in Figure \ref{fig:iters_MP3_nx193_k40_betam1p0_idrs_gmres_F-vs-V}, it can be seen that the convergence of Krylov solvers using the F-cycle is faster than that using the V-cycle. For example, with $np=8$, utilizing the F-cycle yields a wall-clock time of $\SI{336.60}{\second}$ for GMRES and $\SI{101.03}{\second}$ for IDR(4), while implementing the V-cycle results in a wall-clock time of $\SI{796.31}{\second}$
    for GMRES and $\SI{149.16}{\second}$s for IDR(4). In general, the F-cycle involves more computation than the V-cycle due to the additional levels of coarsening and interpolation, but the improved convergence can still make it a viable option in certain cases.


    \subsection{Weak scaling}
    The first parallel property that we are interested in is the weak scalability of our parallel solver. In the experiments, we increase the problem size and the number of processors proportionally to maintain the same grid size for each processor.
    
    The weak scaling results for the 3D closed-off problem are presented in Table \ref{tab: MP1 weak scaling performance}. One can find that parallel GMRES does not show satisfying parallel scalability and consumes more CPU time than IDR(4), especially for large grid sizes. It is because GMRES needs long recurrences and more dot-product operations, hence more global communications. As the problem size increases, the sequential computation load which concerns the process of calculating the Given rotation matrix in GMRES becomes the bottleneck.
    
     For IDR(4), when the number of processors and grid size increases proportionally, that is maintaining the same grid size for each processor, the Wall-clock time remains relatively stable, indicating that our parallel framework can effectively handle large-scale Helmholtz problems with minimized pollution error. Since GMRES does not demonstrate weak scalability for problems requiring a large number of iterations, we will primarily use IDR(4) to solve the following increasingly complex model problems.
    
    \begin{table}[htbp]
        \centering
        \caption{Weak scaling analysis of parallel CSLP-preconditioned GMRES and IDR(4) for 3D Closed-off problem, $k=40$.}
        \scalebox{0.9}{
        \begin{tabular}{r|c|ccc|ccc}
            grid &      & \multicolumn{2}{c}{GMRES}       && \multicolumn{2}{c}{IDR(4)} &\\
            size & $np$ & \#Matvec & $t (s)$              && \#Matvec & $t (s)$  &\\ \hline
            129$\times$129$\times$129 & 8  &146 &30.19    && 199      & 21.03  &\\
            193$\times$193$\times$193 & 27 &137 &80.10    && 187      & 26.93  &\\
            257$\times$257$\times$257 & 64 &143 &97.20    && 190      & 25.22  &\\
            \end{tabular}
        }
        \label{tab: MP1 weak scaling performance}
    \end{table}
    
    
    \subsection{Strong scaling}
    In this section, we perform numerical experiments on a fixed problem size while varying the number of processors. The speedup and parallel efficiency are then calculated to assess the strong scaling of our parallel solver.

    \subsubsection{3D closed-off problem}
    We initiate our analysis by comparing the strong scaling performance of our parallel framework with a PETSc \cite{petsc-user-ref,petsc-efficient} implementation, specifically employing PETSc version 3.19.0 with complex support and optimized mode enabled. The application is compiled and tested in the same environment as our program. We use the CSLP-preconditioned GMRES algorithm to solve the 3D closed-off problem with a grid size of $129\times129\times129$ and wavenumber $k=40$. In the PETSc implementation, matrices are explicitly constructed, including the matrix defining the linear system and the matrix employed for constructing the preconditioner (CSLP). The CSLP is implemented in a so-called Shell preconditioner of PETSc, where CSLP is inverted approximately by a default preconditioned GMRES solver. Both methods exhibit similar convergence properties. The number of iterations required for the PETSc implementation is $155$ for sequential computing and that of our implementation is $146$. Table \ref{tab: MP1 compare petsc} presents the average wall-clock time for each outer GMRES iteration and the corresponding strong scaling. Our implementation demonstrates less time consumption and better parallel performance.

    \begin{table}[htbp]
        \centering
        \caption{Comparisons with PETSc implementation for one GMRES iteration. Parallel CSLP-preconditioned GMRES is used to solve the 3D Closed-off problem with grid size $129\times 129 \times 129$, $k=40$.}
        \scalebox{0.9}{
        \begin{tabular}{c|cccc|cccc}
                  & \multicolumn{3}{c}{PETSc}  &  & \multicolumn{3}{c}{Present} &\\
             $np$ &  $t (s)$  & $S_p$  & $E_p$ &  & $t (s)$  & $S_p$  & $E_p$   &  \\ \hline
               1  &    30.698   &  -     &  -    &  &  1.310   & -      & -       &   \\
               4  &    8.802    &  3.49  & 0.87  &  &  0.334   & 3.92   & 0.97    &   \\
               8  &    5.280    &  5.81  & 0.73  &  &  0.205   & 6.37   & 0.80    &   \\
               16 &    2.759    &  11.13 & 0.69  &  &  0.117   & 11.20  & 0.70    &   \\
            \end{tabular}
        }
        \label{tab: MP1 compare petsc}
    \end{table}

    \subsubsection{3D Wedge problem}
    For the 3D Wedge problem, the performance of the parallel CSLP-preconditioned IDR(4)  on one compute node and multiple compute nodes are shown in Table \ref{tab: MP3 f40 nx193 Sp and Ep} and Table \ref{tab: MP3 f80 nx385 Sp and Ep}, respectively. The results show that the parallel framework exhibits good performance for non-constant wavenumbers and is scalable across multiple compute nodes as well as within a single node. With the growing number of processors, we can observe a moderate decrease in parallel efficiency. This is mainly due to increased communication, which leads to a decrease in the ratio of computation/communication.
    
    \begin{table}[htbp]
        \centering
        \caption{Performance of the parallel CSLP-preconditioned IDR(4) for 3D Wedge problem with grid size $193\times321\times193$ at $f=\SI{40}{\hertz}$.}
        \label{tab: MP3 f40 nx193 Sp and Ep}
        \scalebox{0.89}{
        \begin{tabular}{c|c|c|c|c}
          
          $npx \times npy \times npz$ & \#Matvec & $t (s)$   & $S_p$    & $E_p$   \\ \hline
          1   $\times$  1  $\times$  1  & 395 & 1033.11 &       &      \\
          1   $\times$  2  $\times$  1  & 421 & 557.63  & 1.85  & 0.93 \\
          2   $\times$  2  $\times$  2  & 377 & 152.88  & 6.76  & 0.84 \\
          
        \end{tabular}
      }
    \end{table}
    
    \begin{table}[htp]
        \centering
        \caption{Performance of the parallel CSLP-preconditioned IDR(4) for 3D Wedge problem with grid size $385\times641\times385$ at $f=\SI{80}{\hertz}$.}
        \label{tab: MP3 f80 nx385 Sp and Ep}
        \scalebox{0.89}{
        \begin{tabular}{c|c|c|c|c|c}
          
          $npx \times npy \times npz$& Nodes& \#Matvec & $t (s)$    & $S_p$   & $E_p$   \\ \hline
          4$\times$4$\times$3 & 1 & 835 & 1313.82 &      &      \\
          4$\times$6$\times$4 & 2 & 821 & 952.94  & 1.38 & 0.69 \\
          6$\times$8$\times$4 & 4 & 825 & 418.74  & 3.14 & 0.78 \\
          6$\times$8$\times$6 & 6 & 832 & 298.33  & 4.40 & 0.73 \\
          
        \end{tabular}
      }
    \end{table}

    \subsubsection{3D SEG/EAGE salt model}

    We are interested in the scalability properties of our parallel framework in realistic applications. Hence, we solve this model problem at a fixed frequency on a growing number of compute nodes. Table \ref{tab: MP4 f5 nx641 Sp and Ep DelftBlue} collects the number of matrix-vector multiplications and Wall-clock time versus the number of compute nodes. The fact that solving a linear system with approximately 79.3 million unknowns within hundreds of seconds shows the capability of solving a realistic 3D high-frequency problem with limited memory and time consumption. One can also find the nice property that the number of matrix-vector multiplications keeps independent of the number of compute nodes. Taking the Wall-clock time on a single computing node as a reference, we can observe fairly good scaled parallel efficiency (around $0.8$) for such a large-scale complex model.

    \begin{table}[htbp]
    \centering
    \caption{Performance of the parallel CSLP-preconditioned IDR(4) on DelftBlue for 3D SEG/EAGE Salt Model with grid size $641\times641\times193$ at $f=\SI{5}{\hertz}$.}
    \label{tab: MP4 f5 nx641 Sp and Ep DelftBlue}
    \scalebox{0.89}{
        \begin{tabular}{c|c|c|c|c|c}
           
            $npx \times npy \times npz$ & Nodes & \#Matvec & $t (s)$    & $S_p$   & $E_p$   \\ \hline
            6$\times$4$\times$2 & 1 & 413 & 897.25 &      &      \\
            6$\times$8$\times$2 & 2 & 423 & 510.56 & 1.76 & 0.88 \\
            6$\times$8$\times$4 & 4 & 423 & 298.86 & 3.00 & 0.75 \\
            9$\times$8$\times$4 & 6 & 404 & 203.31 & 4.41 & 0.74 \\
           
        \end{tabular}
    }
    \end{table}
    
    Without loss of generality, we also tried to obtain the parallel performance of our method on a different platform as a complement to this study. In addition to DelftBlue (which includes 48 cores per compute node), we used a commercial supercomputer named Magic Cube 3\footnote{Magic Cube 3: \url{https://www.ssc.net.cn/en/resource-hardware.html}} (which includes 32 cores per compute node) to perform this numerical experiment. Magic Cube 3 is a supercomputer managed by Shanghai Supercomputer Center. It runs on CentOS Linux release 7.5. Each compute node is equipped with two Intel Xeon Gold 6142 processors with 16 cores at 2.6 GHz, and 12$\times$16G DDR4 2666MHz ECC REG memory. The cluster is equipped with an Intel Omni-Path high-speed network with a transmission bandwidth of 100 Gb/s for interconnectivity between compute nodes and storage systems. On Magic Cube 3, Intel Fortran 17.0.4 and Intel MPI 17.4.239 will be used instead.
    
    As shown in Table \ref{tab: MP4 f5 nx641 Sp and Ep Magic Cube}, we can still achieve satisfactory parallel performance. A moderate decrease in terms of the scaled parallel efficiency should be due to the different bandwidths of the platforms. These results indicate that our parallel framework can be used on different computational platforms. This adaptability is crucial for realistic applications.

    \begin{table}[hbtp]
        \centering
        \caption{Performance of the parallel CSLP-preconditioned IDR(4) on DelftBlue for 3D SEG/EAGE Salt Model with grid size $641\times641\times193$ at $f=\SI{5}{\hertz}$.}
        \label{tab: MP4 f5 nx641 Sp and Ep Magic Cube}
        \scalebox{0.89}{
            \begin{tabular}{c|c|c|c|c|c}
                
                $npx \times npy \times npz$ & Nodes & \#Matvec & $t (s)$    & $S_p$   & $E_p$   \\ \hline
                4   $\times$  4  $\times$  2  & 1 & 405      & 505.14  &      &      \\
                4   $\times$  4  $\times$  4  & 2 & 418      & 287.60  & 1.76 & 0.88 \\
                8   $\times$  8  $\times$  2  & 4 & 390      & 155.64  & 3.25 & 0.81 \\
               
            \end{tabular}
        }
    \end{table}

    In summary, our parallel solver exhibits good weak and strong scaling for a variety of test problems, including large-scale, complex, and realistic applications. The results demonstrate the solver's potential for solving challenging 3D large-scale heterogeneous Helmholtz problems with limited memory and time consumption.

\section{Conclusions} \label{sec:concls}
    In this paper, we developed a matrix-free parallel framework of CSLP-preconditioned Krylov subspace methods, such as GMRES, Bi-CGSTAB, and IDR(s), for 3D large-scale Helmholtz problems in heterogeneous media. The preconditioning operator is approximately inverted by a standard 3D multigrid method. We validate the numerical accuracy by comparisons with an analytical result as well as observations of the wave pattern. Both weak and strong scaling properties of our parallel framework for typical non-constant wavenumber model problems are studied. Additionally, a POP performance analysis has been provided in the appendix to address the bottleneck of this framework.

    To sum up, our research presents a novel matrix-free parallelization of the CSLP preconditioner for the Helmholtz problems. Our parallel implementation maintains good convergence properties without excessive memory, making it an effective alternative to traditional matrix-based preconditioners. Moreover, our work provides a robust and scalable matrix-free parallel-computing framework for solving the Helmholtz problem in increasingly complex 3D scenarios, allowing for the use of various Krylov subspace methods and multigrid methods. Its weak scalability makes it possible to solve the Helmholtz problems with minimized pollution error by using a very large grid size. It also provides an implementation direction for researchers to further develop parallel scalable iterative solvers with wavenumber-independent convergence. Finally, we demonstrate the effectiveness and scalability of our solver/preconditioner combinations in solving the Helmholtz equation on a large-scale parallel architecture, which can be beneficial to engineers to solve large-scale heterogeneous Helmholtz problems.

    Future work can focus on further improving the parallel efficiency of the solver and exploring the use of advanced preconditioning techniques to enhance the convergence of Krylov-based iterative methods.
    
    \section*{Acknowledgments}
    We would like to acknowledge the support of the CSC scholarship (No. 202006230087).

    \appendix
    \section{POP Performance Analysis} \label{appendix_POP}
	With traditional performance metrics such as speed-up and efficiency, it is difficult to understand the actual execution behavior of a parallel program and identify the cause of poor performance and where it occurs. In order to explore the bottleneck of our parallel framework and which part can be improved, this section will consider the performance assessment of our code, using the methodology of the Performance Optimisation and Productivity (POP) provided by the EU HPC Centre of Excellence (CoE)\footnote{POP CoE: \url{https://www.pop-coe.eu}}. The POP methodology can help us build a quantitative picture of application behavior by a set of POP performance metrics, including parallel efficiency (PE), load balance (LB), communication efficiency (CommE), serialization efficiency (SerE), transfer efficiency (TE) and so on. The metrics are computed as efficiencies ranging from 0 to 1, where higher values are more desirable. In general, efficiencies $~0.8$ are considered acceptable, while lower values signal performance concerns that warrant further investigation. We use the following open-source tools: Score-P \cite{ScoreP2011} for profiling and tracing, Scalasca \cite{Scalasca2010} for extended analyses, and CUBE \cite{Cube2021} for presentation.

	Table \ref{tab: Wedge nx385 idr4 pop metric} summarizes the performance assessment of the parallel CSLP-preconditioned IDR(4) on different numbers of compute nodes of DelftBlue for the 3D Wedge model problem with grid size $385 \times 641 \times 385$ at $\SI{80}{\hertz}$. The parallel efficiency drops from $0.83$ to $0.5$ when the number of processes reaches $288$. We find that among the two factors contributing to parallel efficiency, a high load balance is maintained, but communication efficiency decreases. Among the two aspects of communication efficiency, excellent transfer efficiency is maintained, while the most significant inefficiency is serialization, which concerns processes waiting at communication points due to temporal imbalance.
	\begin{table}[H]
		\centering
		\caption{Execution efficiency of the parallel CSLP-preconditioned IDR(4) on different numbers of compute nodes of DelftBlue for 3D Wedge model problem with grid size $385 \times 641 \times 385$ at $\SI{80}{\hertz}$}
		\label{tab: Wedge nx385 idr4 pop metric}
			\begin{threeparttable}
			\begin{tabular}{lccc}
				\toprule
				Compute nodes                              & 1    & 2    & 6    \\
				$np$                       & 48   & 96   & 288  \\ \hline
				Parallel Efficiency (PE)\tnote{1}                  & 0.83 & 0.78 & 0.32 \\
				+ Load Balance (LB) \tnote{2}                    & 0.91 & 0.82 & 0.78 \\
				+ Communication Efficiency (CommE)\tnote{3}      & 0.92 & 0.96 & 0.41 \\
				\quad ++ Serialisation Efficiency (SerE) \tnote{4}  & 0.92 & 0.96 & 0.42 \\
				\quad ++ Transfer Efficiency (TE) \tnote{5}        & 0.99 & 1.00 & 0.97 \\ \bottomrule
			\end{tabular}
			\begin{tablenotes}
				\item[1] \textit{Parallel efficiency} is the ratio of mean computation time to total runtime of all processes
				\item[2] \textit{Load balance} is the mean/maximum ratio of computation time outside of MPI
				\item[3] \textit{Communication efficiency} is the ratio of maximum computation time to total runtime.
				\item[4] \textit{Serialisation efficiency} is estimated from idle time within communications where no data is transferred.
				\item[5] \textit{Transfer efficiency} relates to essential time spent in data transfers.
			\end{tablenotes}
		\end{threeparttable}
	\end{table}

    To determine which part contributes to the poor serialization efficiency, we further study the performance of the precondition, matrix-vector multiplication, and dot-product operations, as shown in Tables \ref{tab: Wedge nx385 cslp pop metric}, \ref{tab: Wedge nx385 matvec pop metric}, and \ref{tab: Wedge nx385 dotprod pop metric}, respectively. It can be seen that parallel matrix-vector multiplication can maintain fairly good efficiency. The precondition component, exhibiting similar behavior to the whole framework in Table \ref{tab: Wedge nx385 idr4 pop metric}, can be considered as the main factor that affects the overall efficiency. The results in Table \ref{tab: Wedge nx385 dotprod pop metric} reveal that the dot product operation is one of the main reasons for the low serialization efficiency. Thus, we can conclude that the preconditioning step becomes the bottleneck because it uses full GMRES to solve coarse grid problems, which involve numerous inner product operations. Further efficiency optimization in this direction can be implemented in future work.
	\begin{table}[H]
		\centering
		\caption{Execution efficiency of the precondition part of parallel IDR(4) on different numbers of compute nodes of DelftBlue}
		\label{tab: Wedge nx385 cslp pop metric}
			\begin{tabular}{lccc}
				\toprule
				Compute nodes                              & 1    & 2    & 6    \\
				$np$                       & 48   & 96   & 288  \\ \hline
				Parallel Efficiency (PE)                  & 0.80 & 0.79 & 0.25 \\
				+ Load Balance (LB)                       & 0.90 & 0.86 & 0.79 \\
				+ Communication Efficiency (CommE)        & 0.89 & 0.92 & 0.31 \\
				\quad ++ Serialisation Efficiency (SerE)  & 0.90 & 0.92 & 0.32 \\
				\quad ++ Transfer Efficiency (TE)         & 0.99 & 0.99 & 0.97 \\ \bottomrule
			\end{tabular}
	\end{table}

	\begin{table}[H]
		\centering
		\caption{Execution efficiency of the matrix-vector multiplications of parallel IDR(4) on different numbers of compute nodes of DelftBlue}
		\label{tab: Wedge nx385 matvec pop metric}
			\begin{tabular}{lccc}
				\toprule
				Compute nodes                              & 1    & 2    & 6    \\
				$np$                       & 48   & 96   & 288  \\ \hline
				Parallel Efficiency (PE)                  & 0.70 & 0.62 & 0.55 \\
				+ Load Balance (LB)                       & 0.90 & 0.77 & 0.82 \\
				+ Communication Efficiency (CommE)        & 0.77 & 0.81 & 0.67 \\
				\quad ++ Serialisation Efficiency (SerE)  & 0.78 & 0.82 & 0.67 \\
				\quad ++ Transfer Efficiency (TE)         & 0.99 & 0.99 & 0.99 \\ \bottomrule
			\end{tabular}
	\end{table}

	\begin{table}[H]
		\centering
		\caption{Execution efficiency of dot-product operations of parallel IDR(4) on different numbers of compute nodes of DelftBlue}
		\label{tab: Wedge nx385 dotprod pop metric}
			\begin{tabular}{lccc}
				\toprule
				Compute nodes                              & 1    & 2    & 6    \\
				$np$                       & 48   & 96   & 288  \\ \hline
				Parallel Efficiency (PE)                  & 0.53 & 0.21 & 0.18 \\
				+ Load Balance (LB)                       & 0.70 & 0.55 & 0.65 \\
				+ Communication Efficiency (CommE)        & 0.75 & 0.38 & 0.28 \\
				\quad ++ Serialisation Efficiency (SerE)  & 0.76 & 0.38 & 0.29 \\
				\quad ++ Transfer Efficiency (TE)         & 1.00 & 0.99 & 1.00 \\ \bottomrule
			\end{tabular}
	\end{table}

\bibliographystyle{siam}
\bibliography{Jinqiang_Manuscript_NMLSP2022}
\end{document}